\documentclass{article}
\usepackage[applemac]{inputenc}
\usepackage[T1]{fontenc}
\usepackage{lmodern}
\usepackage[a4paper]{geometry}
\usepackage[english]{babel}
\usepackage{graphicx}
\usepackage{onlyamsmath}
\usepackage{amsfonts}
\usepackage{amssymb}
\usepackage{ntheorem}

\newcommand{\e}{\epsilon}
\newcommand{\se}{\sqrt{\epsilon}}
\newcommand{\ga}{\gamma}

\newcommand{\sga}{\sqrt{\gamma}}
\newcommand{\sgap}{\sqrt{\gamma'}}

\newcommand{\D}{\mathcal{D}}
\newcommand{\E}{\mathbb{E}}
\newcommand{\esp}{\mathcal{H}^\omega}

\newcommand{\ko}[1]{k^2 (\omega_{#1})}
\newcommand{\N}[1]{N (\omega{#1})}

\newcommand{\Bh}[2]{\beta_{#1} (\omega{#2})}

\newcommand{\M}{\widehat{p} (\omega,x,z)}

\newcommand{\Mt}[2]{p^\epsilon _{tr}\left(#1,x,#2\right)}

\newcommand{\p}[1]{\widehat{p}_{#1} (\omega,z)}

\newcommand{\ha}[1]{\widehat{a}_{#1} (\omega,z)}
\newcommand{\hae}[1]{\widehat{a}_{#1}^\epsilon (\omega,z)}

\newcommand{\hb}[1]{\widehat{b}_{#1} (\omega,z)}
\newcommand{\hbe}[1]{\widehat{b}_{#1}^\epsilon (\omega,z)}

\newcommand{\hf}{\widehat{f}(\omega)}
\newcommand{\1}{\textbf{1}}
\newcommand{\dz}{\frac{d}{dz}}

\newcommand{\td}[1]{\tilde{d}_{#1}}

\newtheorem{thm}{Theorem}[section]
\newtheorem{prop}[thm]{Proposition}
\newtheorem{lem}[thm]{Lemma}

\theoremstyle{nonumberplain}
\theorembodyfont{\upshape}

\title{Loss of Resolution for the Time Reversal of Waves in Random Underwater Acoustic Channels}

\author{Christophe Gomez\thanks{Department of Mathematics, Stanford University, Building 380, Sloan Hall
Stanford, California 94305 USA (chgomez@math.stanford.edu). Tel: +(1)650-723-1968. Fax: +(1)650-725-4066.}}

\begin{document}

\maketitle

\begin{abstract}

In this paper we analyze a time-reversal experiment in a random underwater acoustic channel. In this kind of waveguide with semi-infinite cross section a propagating field can be decomposed over three kinds of modes: the propagating modes, the radiating modes and the evanescent modes.  Using an asymptotic analysis based on a separation of scales technique we derive the asymptotic form of the the coupled mode power equation for the propagating modes. This approximation is used to compute the transverse profile of the refocused field and show that random inhomogeneities inside the waveguide deteriorate the spatial refocusing. This result, in an underwater acoustic channel context, is in contradiction with the classical results about time-reversal experiment in other configurations, for which randomness in the propagation medium enhances the refocusing.

\end{abstract}

\begin{flushleft}
\end{flushleft}

\begin{flushleft}
\end{flushleft}

\section*{Introduction}

\begin{figure} \begin{center}
\begin{tabular}{cc}\includegraphics*[scale=0.27]{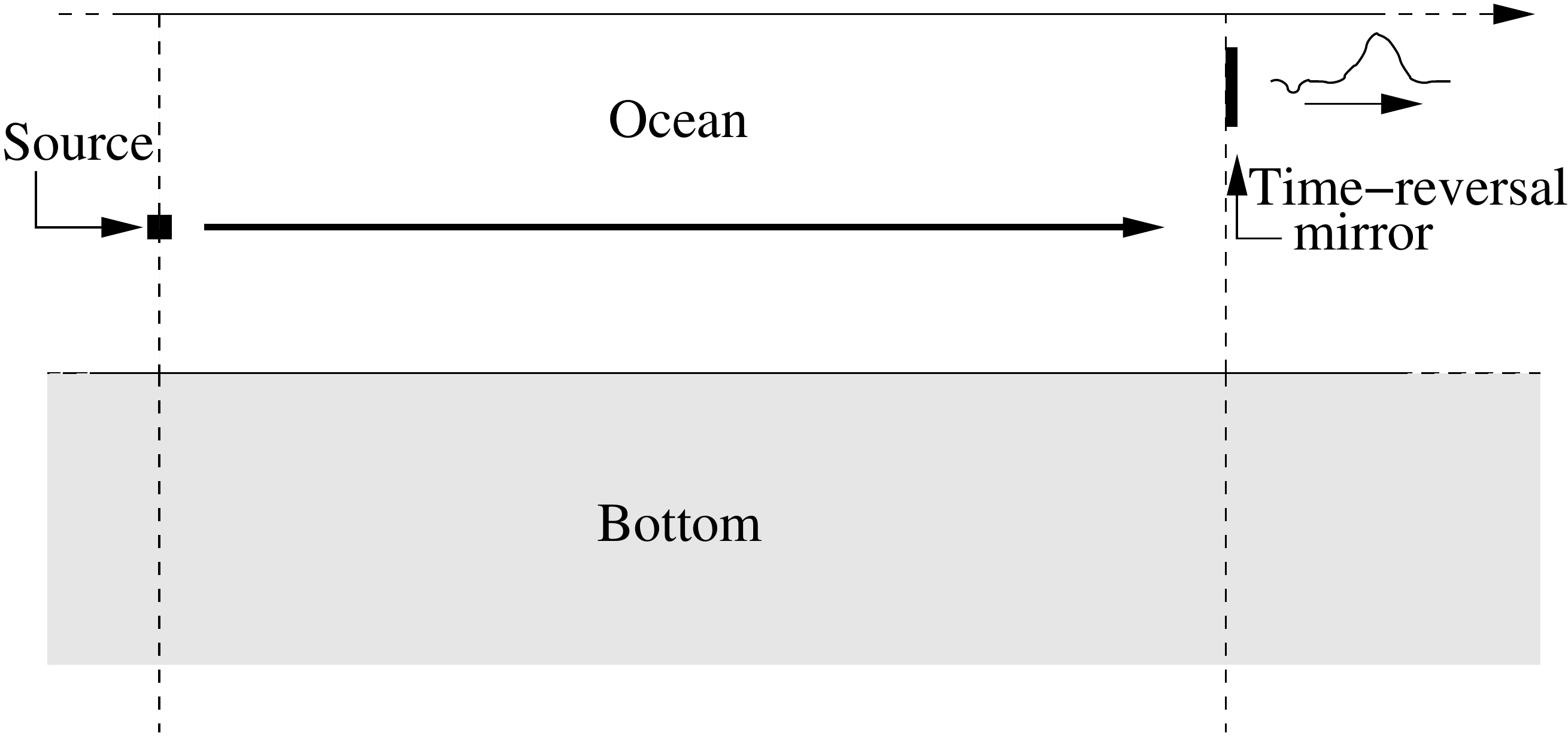}&\includegraphics*[scale=0.27]{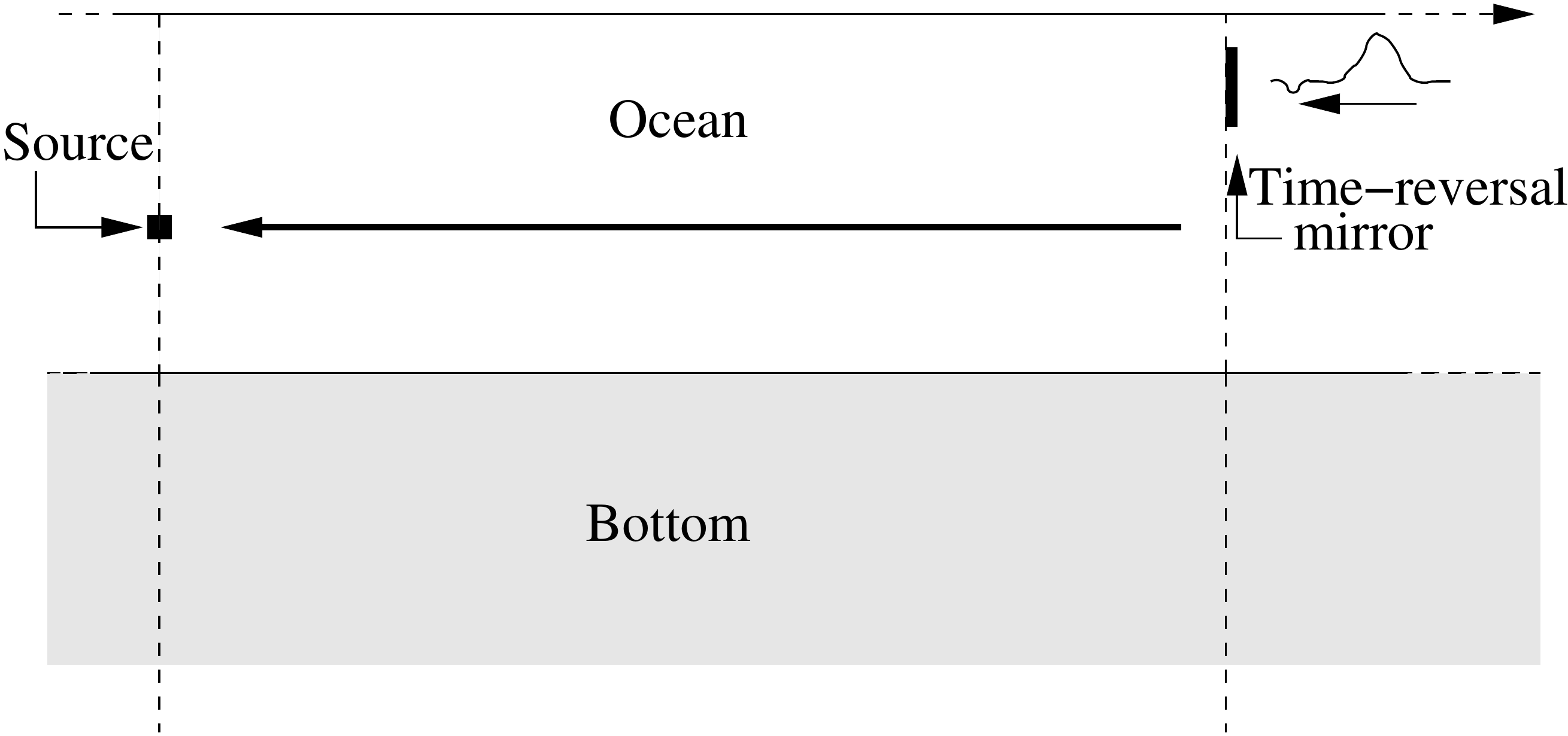}\\
$(a)$ &$(b)$
\end{tabular}
\end{center}
\caption{\label{figtrmintP2p}
Representation of the time-reversal experiment. In $(a)$ we represent the first step of the experiment, and in $(b)$ we represent the second step of the experiment. 
}\end{figure}

The time-reversal experiments of M. Fink and his group in Paris have attracted considerable attention because of the surprising effect of enhanced spatial focusing and time compression in random media. The refocusing properties have numerous interesting applications, in detection, destruction of kidney stones, and wireless communication for instance. This experiment is in two steps. In the first step (see Figure \ref{figtrmintP2p} $(a)$), a source sends a pulse into a medium. The wave propagates and is recorded by a device called a time-reversal mirror. A time-reversal mirror is a device that can receive a signal, record it, and resend it time-reversed into the medium. In the second step (see Figure \ref{figtrmintP2p} $(b)$), the wave emitted by the time-reversal mirror has the property of refocusing near the original source location. However, surprisingly, it has been observed that random inhomogeneities enhance refocusing \cite{fink1,fink2}. 

Time-reversal refocusing in one-dimensional propagation media is carried out in \cite{clouet,book}. In three-dimensional randomly layered media \cite{trfpsource}, in the paraxial approximation \cite{papanicolaou3,papanicolaou2,papanicolaou4}, and in random waveguides \cite{book,papa}, it has been shown that the focal spot can be smaller than the Rayleigh resolution formula $\lambda L/D$ (where $\lambda$ is the carrier wavelength, $L$ is the propagation distance, and $D$ is the mirror diameter), but the focal spot is still larger than the diffraction limit $\lambda/2$. Moreover, in \cite{gomez} the author propose a setup in a waveguide with random perturbations in the vicinity of the source in order to obtain a \emph{superresolution} effect, that is to refocus beyond the diffraction limit with a far-field time-reversal mirror. The setup was inspired from \cite{science} describing this superresolution effect experimentally.

Time reversal of a broadband pulse, in the case of a waveguide with a bounded cross-section and Dirichlet boundary conditions, is carried out in \cite{papa} and \cite[Chapter 20]{book}. In underwater acoustics the waveguide model (see Figure \ref{figtrmintP2p}) has a semi-infinite cross section, and therefore a wave field can be decomposed into three kinds of modes: the propagating modes, the evanescent modes, and the radiating modes. Consequently, waveguides with bounded cross-sections does not take into account radiation losses which can be encountered in underwater acoustics. Moreover, the refocusing enhancement of the time-reversal experiment in random waveguides is closely related to the propagation of the propagating modes power in the random medium. For waveguides with bounded cross-sections the total propagating mode power is conserved and even uniformly distributed after long-range propagations \cite{book,papa}, that is why refocusing enhancement can be observed. In underwater acoustics, with a far-field time-reversal mirror only the propagating modes contribute to the focal spot, and it has been shown \cite{gomez2,papanicolaou} that the presence of radiating modes produces an effective dissipation on the propagating modes, which cannot be fully compensated by time reversal. As a result, in contrast with all the results listed above regarding improvement of the time-reversal focusing thanks to the random perturbations in the propagation medium, we show in this paper that random inhomogeneities in a context of underwater acoustics deteriorate the refocusing property. The deterioration of the refocusing property for the time-reversal experiment has already been observed in \cite{ammari}. However, in there context the authors considered a thermo-viscous wave model to incorporate viscosity effects in wave propagation. In our context, the deterioration of the refocusing is only due to the random inhomogeneities and the geometry of the propagation media. Therefore, it is interesting to understand these effects and to characterize them on the focal spot resulting from the time-reversal experiment. The main result of this paper is the analysis of the radiation losses on the refocused wave in the time-reversal experiment in a context of underwater acoustics. In Theorem \ref{transprofP2p} and Theorem \ref{transprof2P2p}, we show that the radiative losses affect the quality of the time-reversal refocusing in two different ways. First, as expected according to the results obtained in \cite{gomez2}, the amplitude of the refocused wave decays exponentially with the propagation distance. Second, the width of the main focal spot increases and converges to an asymptotic value, which is larger than the diffraction limit $\lambda_{oc}/(2\theta)$ obtained in Proposition \ref{prop17P2p} (where $\lambda_{oc}$ is the carrier wavelength in the ocean section with index of refraction $n_1$, and $\theta=\sqrt{1-1/n_1^2}$).

The organization of this paper is as follows. In Section \ref{sect1P2p} we introduce the underwater waveguide model studied in detail in \cite{wilcox}, and in Section \ref{sect2P2p} we present the mode decomposition associated to this model and we derive the coupled mode equations.  In Section \ref{sect4P2p} we study the time-reversal experiment. We describe in a simple way the refocused transverse profile in terms of the solution of the continuous diffusive model obtained in \cite{gomez2}, and describing the mode-power coupling between the propagating and radiating modes.Thanks to this representation we show that the quality of the time-reversal refocusing is deteriorated by the radiative losses in the ocean bottom.

\section{Waveguide Model}\label{sect1P2p}

\begin{figure}\begin{center}
\includegraphics*[scale=0.4]{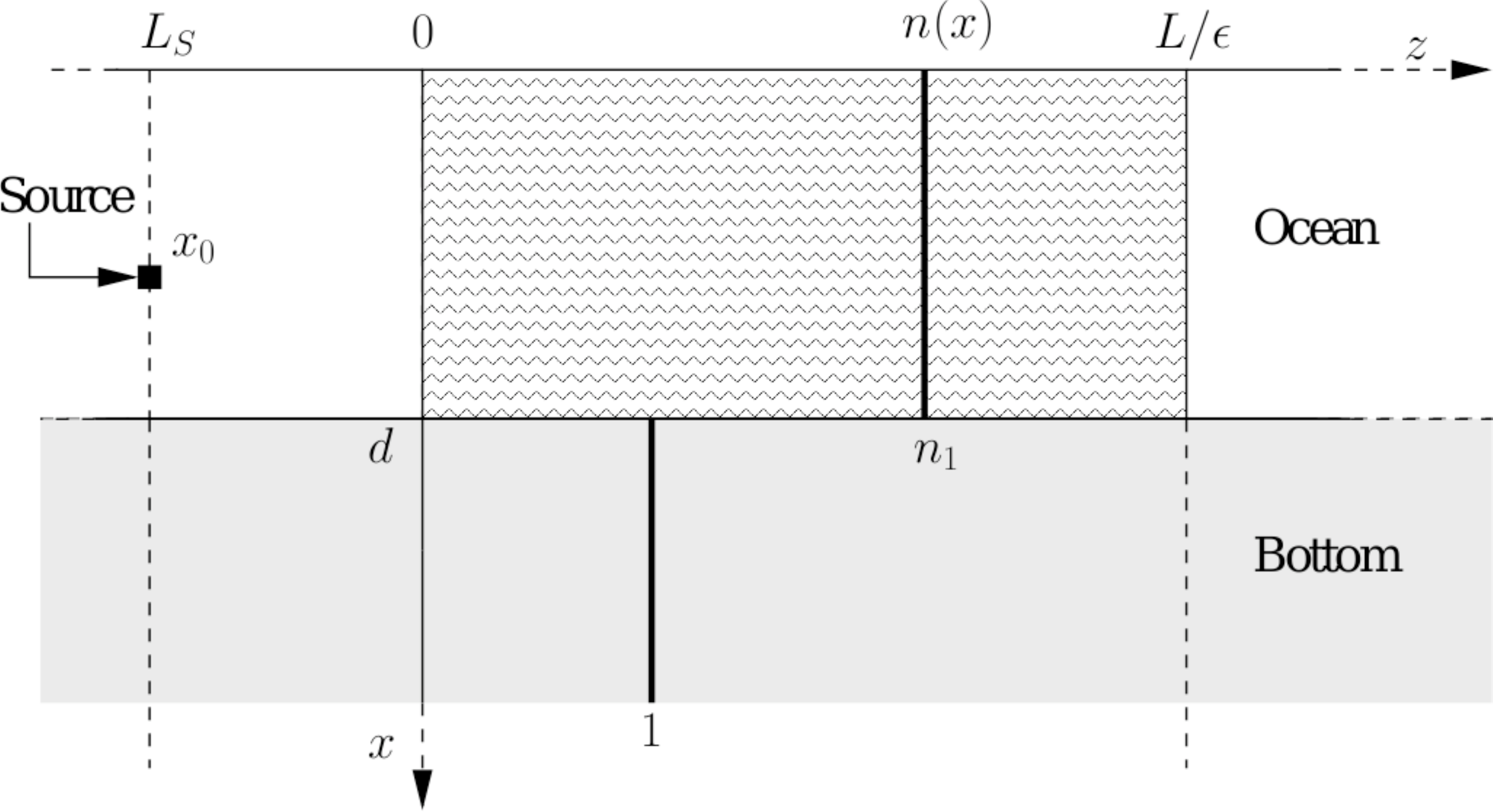}
\end{center}
\caption{\label{fig2}Illustration of the shallow-water waveguide model.}
\end{figure}

We consider a two-dimensional linear acoustic waveguide model. The conservation equations of mass and linear momentum are given by
\begin{equation}\label{conservationP2p}
\begin{split}
\rho ( x,z)\frac{\partial \textbf{u}}{\partial t} + \nabla p &=\textbf{F}^\e_q, \\
\frac{1}{K (x,z)} \frac{\partial p}{\partial t} + \nabla . \textbf{u} &=0, 
\end{split}
\end{equation}
where $p$ is the acoustic pressure, $\textbf{u}$ is the acoustic velocity, $\rho$ is the density of the medium, $K$ is the bulk modulus, and the source is modeled by the forcing term $\textbf{F}^\e_q(t,x,z)$
given by 
\begin{equation*}  \textbf{F}^\e (t,x,z)=\Psi^{\e}(t,x)\delta(z-L_S)\textbf{e}_{z}.\end{equation*} 
The third coordinate $z$ represents the propagation axis along the waveguide. The transverse section of the waveguide is the semi-infinite interval $[0,+\infty)$, and $x \in [0,+\infty)$ represents the transverse coordinate.
Here, $\textbf{F}^\e$ represents a point source localized at $z=L_S$, pointing in the $z$-direction, with temporal and transverse profile given by  $\Psi^\e ( t,x)$. Let $d>0$ be the bottom of the underwater waveguide, the medium parameters are given by
\[\begin{split}
\frac{1}{K(x,z)} & =  \left\{ \begin{array}{ccl} 
                                            \frac{1}{\bar{K}}\left( n^2(x)+\sqrt{\e} V(x,z) \right) & \text{ if }  &  x\in [0,d],\quad z\in [0,L/ \e] \\
                                             \frac{1}{\bar{K}}n^2(x) & \text{ if }  & \left\{\begin{array}{l} x\in[0,+\infty),\,z\in (-\infty,0)\cup(L/\e,+\infty)\\ \text{or}\\ x\in (d,+\infty),\,z\in(-\infty,+\infty). \end{array} \right.
                                          \end{array} \right. \\
\rho(x,z)&=  \bar{\rho}\quad \text{ if }\quad  x\in [0,+\infty),\,z\in\mathbb{R}, \\
\end{split}\] 
and where $V$ is a stochastic process describing the random perturbation of the propagation medium (see Figure \ref{fig2} for an illustration of the random underwater waveguide model). In this paper we consider the Pekeris waveguide model. This kind of model has been studied for half a century \cite{pekeris} and in this model the index of refraction $n(x)$ is given by
\begin{equation*}
n(x)=\left\{\begin{array}{lcl} n_1>1& \text{if}& x\in[0,d)\\
                                                                 1& \text{if}& x\in[d,+\infty).
\end{array}\right.\end{equation*}
The Pekeris profile models an ocean with a constant sound speed, and where $d$ represents the ocean depth. Conditions corresponding to the Pekeris model can be found during the winter in Earth's mid latitudes and in water shallower than about $30$ meters. 

The perturbation $V$ is assumed to be a continuous real-valued zero-mean stationary stochastic process with $\phi$-mixing properties \cite{kushner}. More precisely, let 
\[\mathcal{F}_u=\mathcal{F}_{0,u}=\sigma(V(x,z), \,x\in[0,d],\quad0\leq z\leq u )\quad\text{and}\quad \mathcal{F}_{u,+\infty}=\sigma(V(x,z), \,x\in[0,d],\quad,u\leq z ),\] 
we assume that
\[\sup_{\substack{v\geq 0\\ A\in\mathcal{F}_{u+v,+\infty}\\ B\in\mathcal{F}_{0,v}}}\lvert \mathbb{P}(A\vert B)-\mathbb{P}(A)\rvert \leq \phi(u).\]
This $\phi$-mixing property describes the decorrelating behavior of the random perturbation $V$ through the nonnegative function $\phi\in L^1(\mathbb{R})\cap L^{1/2}(\mathbb{R})$ characterizing the decorrelation speed. Throughout this paper, for the sake of simplicity and for explicit computations in Section \ref{mrfnrlintP2p} and Section \ref{mrfrsintP2p}, we assume that the autocorrelation function of the random perturbation $V$ is given by
\[\E[V(x_1,z_1)V(x_2,z_2)]=\gamma_0(x_1,x_2)e^{-a\lvert z_1-z_2\rvert}\quad \forall (z_1,z_2,x_1,x_2)\in[0,+\infty)^2\times [0,d]^2.\]

From the conservation equations \eqref{conservationP2p}, we derive the wave equation for the pressure field,
\[\Delta p - \frac{1}{ c(x,z) ^2 }\frac{\partial ^2  p }{\partial t^2} = \nabla .\textbf{F}^\e,\]
where $c(x,z)=\sqrt{K(x,z)/\rho(x,z)}$ is the sound speed profile, $\Delta = \partial ^2 _x + \partial ^2 _{z}$, and $c= \sqrt{\bar{K}/\bar{\rho}}$. In underwater acoustics  it is natural to use a pressure-release condition since the density of air is very small compared to the density of water. As a result, the pressure is very weak outside the waveguide, and by continuity, the pressure at the free surface $x=0$ is zero. This consideration leads us to consider the Dirichlet boundary conditions
\begin{equation*} 
p (t,0,z)=0 \quad \forall (t,z)\in [0,+\infty)  \times \mathbb{R}.
\end{equation*}

To study the focusing property of the time-reversal experiment we need to understand the wave propagation in the random medium. To do that, we use a separation of scale technique introduced by G. Papanicolaou and his coauthors in \cite{asch} for instance. The important scale parameters in our problem are: the wavelength, the correlation length and the standard deviation  of the medium inhomogeneities, the propagation distance, and bandwidth of the pulse. This last scale parameter plays a key role in the statistical stability of the time-reversal experiment.

First of all, the model of wave propagation considered in this paper is a linear models, so that the pressure $p(t,x,z)$ can be expressed as the superposition of monochromatic waves by taking its Fourier transform. Here, the Fourier transform and the inverse Fourier transform, with respect to time, are defined by
\[ \hf =\int f(t)e^{i \omega t} dt, \quad f(t)=\frac{1}{2 \pi} \int \hf e^{-i \omega t} d\omega.\] 
As a result, in the half-space $z>L_S$ (resp., $z<L_S$), we get that $\M$ satisfies the time-harmonic wave equation without source term
\begin{equation}\label{helmo2P2} \partial ^2 _z\M+\partial ^2 _x\M +\ko{}n^2 (x)\M+\se V(x,z)\M \1_{[0,d]}(x)\1_{[0,L/\e]}(z)=0, \end{equation}
where $k(\omega)=\omega/c$ is the wavenumber, and with Dirichlet boundary conditions $\widehat{p}(\omega,0,z)=0$ $\forall z$. 
The source term does not appear in \eqref{helmo2P2} but induces the following jump conditions for the pressure field across the plane $z=L_S$
\begin{equation}\label{jumpscondP2}\begin{array}{ccl}
\widehat{p}(\omega,x,L_S^+)-\widehat{p}(\omega,x,L_S^-)&=& \widehat{\Psi}(\omega,x),\\
\partial_z \widehat{p}(\omega,x,L_S^+)-\partial_z \widehat{p}(\omega,x,L_S^-)&=&0.
\end{array}
\end{equation}
To study \eqref{helmo2P2}, we consider this equation as an operational differential equation
 \[
 \frac{d^2}{d z^2} \widehat{p}(\omega,.,z)+R(\omega) \big(\widehat{p}(\omega,.,z)\big)+\se V(\cdot,z)\widehat{p}(\omega,.,z) \1_{[0,d]}(\cdot)\1_{[0,L/\e]}(z)=0 \]
in $H=L^2(0,+\infty)$, where $R(\omega)$ is an unbounded operator on $H$ with domain
\[ \mathcal{D}(R(\omega))= H^1 _0 (0,+\infty)\cap H^2(0,+\infty),\]
and defined by
\begin{equation}\label{pekerisop}R(\omega)(y)=\frac{d^2}{ dx^2} y +\ko{}n^2 (x)y\quad \forall y\in \mathcal{D}(R(\omega)).\end{equation}
In the next section we introduce the spectral decomposition of the operator $R(\omega)$ \cite{wilcox}. This decomposition we will be used in what follows to decompose the field $\M$ and then understand the stochastic effects undergoes during the propagation.

\subsection{Spectral Decomposition in Unperturbed Waveguides}\label{spectralP2}

The spectral analysis of the self-adjoint Pekeris operator \eqref{pekerisop} is carried out in \cite{wilcox}. To use this spectral decomposition, we are interested in solutions of \eqref{helmo2P2} such that
\begin{equation*}\begin{split}
\widehat{p}(\omega,.,.)\textbf{1}_{(L_S,+\infty)}(z) &\in \mathcal{C}^0 \Big((L_S,+\infty), H^1 _0 (0,+\infty)\cap H^2(0,+\infty) \Big)\cap \mathcal{C}^2 \Big((L_S,+\infty),H\Big),\\
\widehat{p}(\omega,.,.)\textbf{1}_{(-\infty,L_S)}(z)&\in \mathcal{C}^0 \Big((-\infty,L_S), H^1 _0 (0,+\infty)\cap H^2(0,+\infty) \Big)\cap \mathcal{C}^2 \Big((-\infty,L_S),H\Big).
\end{split}\end{equation*}  
According to \cite{wilcox}, the spectrum of the unbounded operator \eqref{pekerisop} is given by
\[
Sp\big(R(\omega)\big)=\left(-\infty,\ko{}\right]\cup\big\{ \beta^2 _{\N{}}(\omega),\dots,\beta_1^2 (\omega)\big\}.\]
The continuous part of the spectrum comes from the fact that our waveguide model is semi-infinite. For the discrete part, the modal wavenumber $\beta_j (\omega)$ are positive and
\[\ko{}<\beta^2 _{\N{}}(\omega)<\cdots<\beta_1^2 (\omega)<n_1^2 \ko{}.\]
Regarding the spectral decomposition, there exists a resolution of the identity $\Pi_\omega$ of $R(\omega)$ such that $\forall y\in H$, $\forall r\in\mathbb{R}$,

\begin{equation*}\begin{split}
 \Pi_\omega(r,+\infty)(y) (x) =&  \sum_{j=1}^{\N{}} \big<y,\phi_j(\omega,.)\big>_H \phi_j(\omega,x)\textbf{1}_{(r,+\infty)}\left(\Bh{j}{}^2\right) \\
                                                                      & + \int_{r}^{\ko{}}  \big<y,\phi_\ga(\omega,.) \big>_H \phi_\ga (\omega,x)d\ga \textbf{1}_{\left(-\infty,\ko{} \right)}(r),
  \end{split}\end{equation*}
and $\forall y\in \mathcal{D}(R(\omega))$, $\forall r\in\mathbb{R}$,  
  \begin{equation*}\begin{split}                                                                    
 \Pi_\omega(r,+\infty)(R(\omega)(y)) (x) =&  \sum_{j=1}^{\N{}}\Bh{j}{}^2 \big<y,\phi_j(\omega,.)\big>_H \phi_j(\omega,x)\textbf{1}_{(r,+\infty)}\left(\Bh{j}{}^2\right) \\
                                                                      & + \int_{r}^{\ko{}}  \ga \big<y,\phi_\ga(\omega,x) \big>_H \phi_\ga (\omega,x)d\ga \textbf{1}_{\left(-\infty,\ko{} \right)}(r).
\end{split}\end{equation*}   
Let us describe more closely the discrete and the continuous part of the decompositions.

\paragraph{Discrete part of the decomposition}

$\forall j\in \big\{1,\dots,\N{}\big\}$, the $j$th eigenvector is given in \cite{wilcox} by
\[\phi_j(\omega, x)=\left\{ \begin{array}{ccl}
A_j(\omega)\sin(\sigma_j(\omega) x/d) & \mbox{ if } & 0\leq x \leq d \\
A_j(\omega)\sin(\sigma_j(\omega)  )e^{-\zeta_j (\omega) \frac{x-d}{d}}& \mbox{ if } & d\leq x,  \end{array} \right.\]
where
\[\sigma_j (\omega)=d \sqrt{n_1 ^2 \ko{}-\beta^2 _j( \omega)}, \quad \zeta_j(\omega) =d \sqrt{\beta^2_j(\omega)-\ko{}},\]
and 
\begin{equation}\label{coefajP2}
A_j(\omega)=\sqrt{\frac{2/d}{1+\frac{\sin^2 (\sigma_j (\omega))}{\zeta_j (\omega) }-\frac{\sin(2 \sigma_j (\omega)  )}{2\sigma_j (\omega) } }}.\end{equation}
Here,  $\sigma_1 (\omega),\dots,\sigma_{\N{}}(\omega)$ are the solutions on $(0,n_1k(\omega)d\theta)$ of the following equation, 
\begin{equation}\label{eqvpP2}
\tan(y )=-\frac{y }{\sqrt{( n_1 k d  \theta) ^2-y^2}},
\end{equation}
and such that $0<\sigma_1 (\omega)<\cdots<\sigma_{\N{}}(\omega)< n_1 k(\omega)d \theta$, with  $\theta =\sqrt{1-1/n_1 ^2}$. This last equation admits exactly one solution over each interval of the form $\big(\pi/2+(j-1)\pi, \pi/2 +j\pi \big)$ for $ j \in \{1,\dots, N(\omega) \}$, where 
\[\N{}=\left[ \frac{n_1  k(\omega)d}{\pi} \theta \right],\] 
and $[\cdot]$ stands for the integer part. From \eqref{eqvpP2}, we have the following results \cite{these} which are used to obtain the main result of this paper in Section \ref{mrfrsintP2p}.
\begin{lem}\label{coefgP2}
Let $\alpha\in(1/3,1)$, we have as $\N{} \to +\infty$
\begin{equation*} 
\sup_{j\in \{1,\dots, \N{}-[\N{}^\alpha] -1 \}} \left\lvert \sigma_{j+1}(\omega)-\sigma_j(\omega) -\pi\right\rvert =\mathcal{O}\left( \N{}^{\frac{1}{2}-\frac{3}{2}\alpha } \right).
\end{equation*}

\begin{equation*} 
\sup_{j\in \{1,\dots,\N{}-[\N{}^\alpha ]-2 \}} \left\lvert \sigma_{j+2}(\omega)-2\sigma_{j+1}(\omega)+\sigma_{j}(\omega)\big)\right\rvert =\mathcal{O}\left( \N{}^{1-3\alpha } \right).
\end{equation*}
\end{lem}   
Let us note that, $\forall \eta \in[0,1[$, we have
\begin{equation}\label{approxvp}\sup_{j\in\{1,\dots,\N{}^{\alpha}\}}\lvert \sigma_j(\omega)- j\pi \rvert=\mathcal{O}(\N{}^{\alpha-1}),\end{equation}
and
\[\lim_{N(\omega)\to+\infty}\sup_{j\in\{1,\dots,\N{}^{\alpha}\}}\| \phi_j(\omega,\cdot)-\phi_j(\cdot)\|_{H} =0\]
with
\[ \phi_j(x)=\left\{\begin{array}{ccl}
                                            \sqrt{\frac{2}{d}}\sin(j\frac{\pi}{d}x) &\text{ if }& x\in[0,d]\\
                                            0 &\text{ if }& x\geq d.\end{array}\right.  \] 
This result means that in the limit of large number of propagating modes the low order propagating modes are very similar in shape to those of a perfect bounded waveguide with pressure-release boundary conditions at $x=0$ and $x=d$. This approximation does not hold anymore for high order propagating modes, but the results of Lemma \ref{coefgP2} mean that the distribution of solutions of \eqref{eqvpP2} is closed to the distribution of the eigenvalues of the transverse Laplacian associated to a perfect bounded waveguide with pressure-release boundary conditions at $x=0$ and $x=d$.

\paragraph{Continuous part of the decomposition}

For $\ga \in (-\infty, \ko{})$, we have \cite{wilcox}
\[
\begin{split}
\phi_\ga& (\omega, x)=\\ 
&\left\{ 
\begin{array}{ccl}
A_\gamma(\omega) \sin(\eta(\omega)  x/d ) & \mbox{ if } & 0\leq x \leq d \\
A_\gamma(\omega) \left(\sin(\eta(\omega)  )\cos\big(\xi(\omega) \frac{x-d}{d}\big)+\frac{\eta(\omega) }{\xi(\omega) }\cos(\eta(\omega) )\sin\big(\xi(\omega) \frac{x-d}{d}\big)\right)& \mbox{ if } & d\leq x, \end{array}
 \right.
\end{split}
\]
where
\[\eta (\omega) =d\sqrt{n_1 ^2 \ko{}-\gamma }, \quad \xi(\omega)  =d\sqrt{\ko{}-\gamma },\]
and
\[A_\gamma(\omega) =\sqrt{\frac{d \xi(\omega) }{\pi\big(\xi ^2(\omega) \sin^{2}(\eta(\omega))+\eta ^2(\omega) \cos^2 (\eta(\omega))\big)}}.\]
Let us note that $\phi_\ga(\omega, .)$ does not belong to $H$ so that $\big<y,\phi_\ga (\omega,.)\big>_H$ is not defined in the classical way. In fact, we have
\[ \big<y,\phi_\ga (\omega,.) \big>_H=\lim_{M\to+\infty}\int_{0}^M y(x) \phi_\ga (\omega, x)dx\]  
where the limit holds on $L^2\big(-\infty,\ko{}\big)$.
Moreover, according to the following Plancherel equality 
\[ \|y \|^2 _H=\|\Pi_\omega(-\infty,+\infty)(y) \|^2 _H = \sum_{j=1}^{\N{}}\big \lvert \big<y,\phi_j(\omega,.) \big>_H \big \rvert^2+\int_{-\infty}^{\ko{}}\big\lvert \big<y,\phi_\ga (\omega,.) \big>_H\big\rvert^2 d\ga,\]
the map which assigns to every element of $H$ the coefficients of its spectral decomposition
\[
\begin{array}{rcc}
\Theta_\omega:H&\longrightarrow &\mathcal{H}^\omega\\
y& \longrightarrow & \Big(\big(\big<y,\phi_j(\omega,.)\big>_H\big)_{ j=1,...,\N{}},\big(\big<y,\phi_\ga (\omega,.)\big>_H\big)_{\ga\in(-\infty,\ko{})}\Big)
\end{array}
\]
is an isometry, from $H$ onto $\esp=\mathbb{C}^{\N{}}\times L^2\big(-\infty,\ko{}\big)$.

\section{Mode Coupling in Random Waveguides} \label{sect2P2p}

Before describing the time-reversal experiment in our randomly perturbed waveguide model, we need to understand how the wave is perturbed during the propagation through the medium. In this section, we study the random effects produced on the modal decomposition of $\M$ propagating in the perturbed section $[0,L/\e]$. 

Using the resolution of the identity $\Pi_\omega$ associated to Pekeris operator $R(\omega)$, we have
\[\M=\sum_{j=1}^{\N{}} \widehat{p}_j (\omega,z )\phi_j(\omega, x)+\int_{-\infty}^{\ko{}}\widehat{p}_\ga (\omega,z)\phi_\ga (\omega,x)d \ga,\]
where $\widehat{p}(\omega,z)=\Theta _\omega(\widehat{p}(\omega,.,z))$ and $\Theta_\omega$ is defined in Section \ref{spectralP2}. 

For the sake of simplicity in the presentation of the forthcoming asymptotic analysis, we will restrict ourself to solutions of the form 
\begin{equation}\label{decfieldsimpP2p}\M=\sum_{j=1}^{\N{}} \widehat{p}_j (\omega,z )\phi_j(\omega, x)+\int_\xi^{\ko{}}\widehat{p}_\ga (\omega,z)\phi_\ga (\omega,x)d \ga.\end{equation}
This assumption is tantamount to neglecting the role played by the evanescent modes during the propagation in the random medium. Nevertheless, as it has been observed in \cite{gomez2, book} that these modes play no role in the refocusing process. The reason is that these modes only imply a mode-dependent and a frequency-dependent phase modulations without remove any energy from the propagating and radiating modes, However, the dispersion phenomena are compensated by the time-reversal mechanism. Moreover, we assume that $\e\ll\xi$ and therefore we have two distinct scales. We will consider in a first time the asymptotic $\e$ goes to $0$ and in a second time the asymptotic $\xi$ goes to $0$.

\subsection{Coupled Mode Equations}

According to the pressure field decomposition \eqref{decfieldsimpP2p}, we give in this section the coupled mode equations, which describes the coupling mechanism between the amplitudes of the two kinds of modes, propagating and radiating modes. In the random section $[0,L/\e]$, $\widehat{p}(\omega,z)$ satisfies the following coupled equation in $\esp_\xi=\mathbb{C}^{\N{}}\times L^2(\xi,\ko{})$.
\begin{equation}\label{eqdiff2P2p}\begin{split}
\frac{d^2}{dz^2}\p{j}+\beta^2 _j(\omega) \p{j}&+\se\ko{}\sum_{l=1}^{\N{}}C^\omega_{jl}(z)\p{l}\\
            &+\se\ko{} \int_\xi^{\ko{}}C^\omega_{j\ga'}(z)\p{\ga'}d\ga'=0,\\
\frac{d^2}{dz^2}\p{\ga}+\ga \,\,\p{\ga}&+\se\ko{}\sum_{l=1}^{\N{}}C^\omega_{\ga l}(z)\p{l}\\
            &+\se\ko{} \int_\xi^{\ko{}} C^\omega_{\ga \ga'}(z)\p{\ga'}d\ga'=0,            
\end{split}\end{equation}
where the coupling coefficients $C^\omega(z)$ are defined by:
\begin{equation}\label{coefcouplVP2}\begin{split}
C^\omega_{jl}(z)&=\big<\phi_j(\omega,.),\phi_l(\omega,.) V(.,z)\big>_{H}=\int_0 ^d\phi_{j}(\omega, x) \phi_l(\omega, x) V(x,z) dx,\\
C^\omega_{j\ga}(z)&=C_{\ga j}(z)=\big<\phi_j(\omega,.),\phi_\ga(\omega,.) V(.,z)\big>_{H} =\int_0 ^d\phi_{j}(\omega, x) \phi_\ga(\omega, x) V(x,z) dx,\\
C^\omega_{\ga \ga'}(z)&=\big<\phi_\ga(\omega,.),\phi_{\ga'}(\omega,.) V(.,z)\big>_{H}=\int_0 ^d\phi_{\ga}(\omega, x) \phi_{\ga'}(\omega, x) V(x,z) dx.
\end{split}\end{equation}
 
Next, we decompose the wave field $\widehat{p}(\omega,z)$ using the amplitudes of the  generalized right- and left-going modes 
$\widehat{a}(\omega,z)$ and $\widehat{b}(\omega,z)$, which are given by
\begin{equation*}\begin{split}
\p{j}&=\frac{1}{\sqrt{\Bh{j}{}}}\Big( \ha{j} e^{i\Bh{j}{}z} +\hb{j}e^{-i\Bh{j}{}z} \Big), \\
\frac{ d }{dz}\p{j} &= i \sqrt{\Bh{j}{}} \Big( \ha{j} e^{i\Bh{j}{}z} - \hb{j}e^{-i\Bh{j}{}z} \Big),\\
\p{\ga}&= \frac{1}{\ga^{1/4}}\Big( \ha{\ga} e^{i\sga z} +\hb{\ga}e^{-i\sga z} \Big),\\
\frac{ d }{dz}\p{\ga} &= i \ga^{1/4} \Big( \ha{\ga} e^{i\sga z} - \hb{\ga}e^{-i\sga z} \Big)
\end{split}\end{equation*}
$\forall j \in \big\{1,\dots,\N{}\big\}$ and almost every $\ga \in (\xi,\ko{})$. 
From \eqref{eqdiff2P2p}, this decomposition allows us to obtain a first order differential system instead of a second order one, so that we obtain the coupled mode equation in $\esp_\xi \times \esp_\xi$ for the amplitudes $(\widehat{a},\widehat{b})$,
\begin{equation}\label{Eq1}\begin{split} 
\dz \widehat{a}(\omega, z)&=\se\, \textbf{H}^{aa}(\omega,z)\big(\widehat{a}(\omega, z)\big)+\se\, \textbf{H}^{ab}(\omega,z)\big(\widehat{b}(\omega, z)\big)
\end{split}\end{equation} 
\begin{equation}\label{Eq2}\begin{split}
\dz \widehat{b}(\omega, z)&=\se\, \textbf{H}^{ba}(\omega,z)\big(\widehat{a}(\omega, z)\big)+\se\, \textbf{H}^{bb}(\omega,z)\big(\widehat{b}(\omega, z)\big).
\end{split}\end{equation}
This system is complemented with the boundary conditions
\[ \widehat{a}(\omega,0)=\widehat{a}^\e_{0}(\omega)\quad \text{ and } \quad\widehat{b}\left(\omega,\frac{L}{\e}\right)=0\]
where
\begin{equation}\label{initfieldP2}\begin{split}
\widehat{a}^\e_{j,0} (\omega) &= \frac{\sqrt{\Bh{j}{}}}{2} \big<\widehat{\Psi}^\e (\omega,\cdot),\phi_j(\omega)\big>_H e^{-i\Bh{j}{}L_S},\quad\forall j \in \big\{1,\dots,\N{} \big\},\\
\widehat{a}^\e_{\ga,0} (\omega)&=\frac{\ga ^{1/4}}{2} \big<\widehat{\Psi}^\e (\omega,\cdot),\phi_\ga(\omega)\big>_H e^{-i \sga L_S}, \quad\text{for almost every } \ga \in(\xi, \ko{}).
\end{split}\end{equation}
For $j\in \big\{1,\dots,\N{}\big\}$, $\widehat{a}_{j,0}(\omega_0)$ represents the initial amplitude of the $j$th propagating mode, and for $\ga \in(\xi, \ko{})$, $\widehat{a}_{\ga,0}(\omega)$ represents the initial amplitude of the $\ga$th radiating mode at $z=0$. The initial conditions for the right-going mode $\widehat{a}_0(\omega_0)$ comes from \eqref{jumpscondP2}, and the second condition for the left-going modes means that no wave is coming  from the right homogeneous waveguide. The coupling operator $\textbf{H}^{aa}(\omega,z)$, $\textbf{H}^{ab}(\omega,z)$, $\textbf{H}^{ba}(\omega,z)$, and $\textbf{H}^{bb}(\omega,z)$ in \eqref{Eq1} and \eqref{Eq2} are defined by:
\[\begin{split}
\textbf{H}^{aa}_{j}(\omega,z)(y)=\overline{\textbf{H}^{bb}_{j}(\omega, z)}(y)&=\frac{i\ko{}}{2}\Big[\sum_{l=1}^{\N{}} \frac{C^\omega_{jl}(z)}{\sqrt{\beta_j(\omega) \beta_l(\omega)}}y_l e^{i(\beta_l(\omega) -\beta_j(\omega))z}\\
&+ \int_{\xi}^{\ko{}}\frac{C^\omega_{j\ga'}(z)}{\sqrt{\beta_j(\omega) \sgap}}y_{\ga'}e^{i(\sgap -\beta_j(\omega))z}d\ga' \Big],
\end{split}\]
\[\begin{split}
\textbf{H}^{aa}_{\ga}(\omega, z)(y)=\overline{\textbf{H}^{bb}_{\ga}(\omega, z)}(y)&=\frac{i\ko{}}{2}\Big[\sum_{l=1}^{\N{}} \frac{C^\omega_{\ga l}(z)}{\sqrt{\sga \beta_l(\omega)}}y_l e^{i(\beta_l(\omega) -\sga)z}\\
&+ \int_{\xi}^{\ko{}}\frac{C^\omega_{\ga \ga'}(z)}{\ga^{1/4}{\ga'}^{1/4}}y_{\ga'}e^{i(\sgap -\sga)z}d\ga' \Big],
\end{split}\]
\[\begin{split}
\textbf{H}^{ab}_{j}(\omega,z)(y)=\overline{\textbf{H}^{ba}_{j}(\omega, z)}(y)&=\frac{i\ko{}}{2}\Big[\sum_{l=1}^{\N{}} \frac{C^\omega_{jl}(z)}{\sqrt{\beta_j(\omega) \beta_l(\omega)}}y_l e^{-i(\beta_l(\omega) +\beta_j(\omega))z}\\
&+ \int_{\xi}^{\ko{}}\frac{C^\omega_{j\ga'}(z)}{\sqrt{\beta_j(\omega) \sgap}}y_{\ga'}e^{-i(\sgap +\beta_j(\omega))z}d\ga' \Big],
\end{split}\]
\[\begin{split}
\textbf{H}^{ab}_{\ga}(\omega,z)(y)=\overline{\textbf{H}^{ba}_{\ga}(\omega,z)}(y)&=\frac{i\ko{}}{2}\Big[\sum_{l=1}^{\N{}} \frac{C^\omega_{\ga l}(z)}{\sqrt{\sga \beta_l(\omega)}}y_l e^{-i(\beta_l(\omega) + \sga)z}\\
&+ \int_{\xi}^{\ko{}}\frac{C_{\ga \ga'}(z)}{\ga^{1/4}{\ga'}^{1/4}}y_{\ga'}e^{-i(\sgap +\sga)z}d\ga' \Big].
\end{split}\]

 Let us remark that we have the following global conservation relations
\[\begin{split}
\|\widehat{a}(\omega,z)\|^2_{\esp_\xi}-\|\widehat{b}(\omega,z)\|^2_{\esp_\xi}&=\|\widehat{a}(\omega,0)\|^2_{\esp_\xi}-\|\widehat{b}(\omega,0)\|^2_{\esp_\xi}\quad \forall z\in\left[0,L/\e\right],\\
\|\widehat{a}\left(\omega,L/\e\right)\|^2_{\esp_\xi}+\|\widehat{b}(\omega,0)\|^2_{\esp_\xi}&=\|\widehat{a}(\omega,0)\|^2_{\esp_\xi}.
\end{split}\]
However,  in our context, the coupling mechanism between the right- and the left-going mode is not very convenient to study the time-reversal experiment. The asymptotic behavior of the whole coupling mechanism between  the right- and the left-going in a random waveguide with a bounded cross section is carried out in \cite{garnier2}, but this study leads to technical difficulties because of the waveguide geometry in our context. Consequently, for the sake of simplicity we introduce in the following section the forward scattering approximation, which allows us to neglect  the coupling mechanism between  the right- and the left-going under certain conditions.

\subsection{Propagator and Forward Scattering Approximation}

Before introducing this approximation, let us define the rescaled processes according to the size of the random section $[0,L/ \e]$, 
\[  \widehat{a}^\e(\omega, z)= \widehat{a}\Big(\omega, \frac{z}{\e}\Big)\quad\text{and}\quad \widehat{b}^\e(\omega, z)= \widehat{b}\Big(\omega, \frac{z}{\e}\Big) \]
which satisfy in $\esp_\xi$ the rescaled coupled mode equation
\[\begin{split} 
\dz \widehat{a}^\e(\omega, z)&=\frac{1}{\se}\, \textbf{H}^{aa}\left(\omega,\frac{z}{\e}\right)\big(\widehat{a}^\e(\omega, z)\big)+\frac{1}{\se}\, \textbf{H}^{ab}\left(\omega,\frac{z}{\e}\right)\big(\widehat{b}^\e(\omega, z)\big)\\
\dz \widehat{b}^\e(\omega, z)&=\frac{1}{\se}\, \textbf{H}^{ba}\left(\omega,\frac{z}{\e}\right)\big(\widehat{a}^\e(\omega, z)\big)+\frac{1}{\se}\, \textbf{H}^{bb}\left(\omega,\frac{z}{\e}\right)\big(\widehat{b}^\e(\omega, z)\big),
\end{split}\]
with the two-point boundary conditions
\[  \widehat{a}^\e(\omega,0)=\widehat{a}^\e_{0}(\omega)\quad \text{ and } \quad\widehat{b}^\e(\omega, L)=0.\]

The propagator $\textbf{P}^\e (\omega, z)$ is defined has being the unique solution of the following differential equation
\[\dz \textbf{P} (\omega, z) = \frac{1}{\sqrt{\e}} \textbf{H} \left(\omega, \frac{z}{\e} \right)\textbf{P}^\e (\omega, z)\quad  \text{with} \quad \textbf{P}^\e (\omega,0)=Id,\]
so that
\[\begin{bmatrix} \hae{} \\ \hbe{} \end{bmatrix}=\textbf{P}^\e (\omega, z) \begin{bmatrix} \widehat{a}^\e (\omega, 0) \\ \widehat{b}^\e  (\omega, 0) \end{bmatrix}.\]
According to the symmetry of $\textbf{H} (\omega, z)$ the propagator has the following particular form 
\[ \textbf{P}^\e (\omega, z)=\begin{bmatrix}  \textbf{P}^a _\e (\omega, z) & \textbf{P}^b _\e (\omega, z) \\ \overline{\textbf{P}^b _\e (\omega, z)} & \overline{\textbf{P}^a _\e (\omega, z)} \end{bmatrix}.\]
where, $\textbf{P}^a _\e (\omega, z)$ and $\textbf{P}^b _\e (\omega, z)$ are two operators representing respectively the coupling between the right-going modes and the coupling between the right-going and left-going modes.

The forward scattering approximation is widely used in the literature. In this approximation the coupling between forward- and backward-propagating modes is assumed to be negligible compared to the coupling between the forward-propagating modes. The physical explanation of this approximation is as follows. The coupling between a right-going propagating mode and a left-going propagating mode involves coefficients of the form
\[ \int_{0}^{+\infty}\mathbb{E}[C^\omega_{jl}(0)C^\omega_{jl}(z)]\cos\big((\Bh{l}{}+\Bh{j}{})z\big)dz,\] 
where the coefficients $C^\omega(z)$ are defined by \eqref{coefcouplVP2}, and the coupling between two right-going propagating modes or two left-going propagating modes involves coefficients of the form
\[ \int_{0}^{+\infty}\mathbb{E}[C^\omega_{jl}(0)C^\omega_{jl}(z)]\cos\big((\Bh{l}{}-\Bh{j}{})z\big)dz\]
$\forall (j,l) \in \big\{1,\dots,\N{}\big\}^2$. The forward scattering approximation consists in assuming that 
\[ \int_{0}^{+\infty}\mathbb{E}[C^\omega_{jl}(0)C^\omega_{jl}(z)]\cos\big((\Bh{l}{}+\Bh{j}{})z\big)dz=0 \quad \forall (j,l) \in \big\{1,\dots,\N{}\big\}^2,\]
so that $\textbf{P}^b _\e (\omega, z)=0$, that is there is no coupling between right-going and left-going propagating modes. Therefore, this approximation holds if the power spectral density of the process $V$, i.e. the Fourier transform of its $z$-autocorrelation function, possesses a cut-off wavenumber. We refer to \cite{papa, gomez2} for justifications on the validity of this approximation. 
As a result, under this approximation we can neglect the left-going propagating modes in the asymptotic $\e\rightarrow0$, and then consider only the simplified coupled amplitude equation on $[0,L]$
\[\dz \hae{}=\frac{1}{\sqrt{\e}}\textbf{H}^{aa}\left(\omega,\frac{z}{\e}\right)
(\hae{})\quad \text{with} \quad \widehat{a}^\e (\omega,0)=\widehat{a} ^\e_0(\omega).\]
Finally, we introduce the transfer operator $\textbf{T}^{\xi,\e} (\omega,z)$, which is the solution of
\begin{equation}\label{transferP2p}
\dz \textbf{T}^{\xi,\e} (\omega,z)=\frac{1}{\sqrt{\e}} \textbf{H}^{aa} \left(\omega, \frac{z}{\e}\right)
\textbf{T}^{\xi,\e} (\omega,z)\quad \text{with} \quad \textbf{T}^{\xi,\e} (\omega,0)= Id.
\end{equation}
From this equation, one can easily check that the transfer operator $\textbf{T}^{\xi,\e} (\omega,z)$ is unitary since $\textbf{H}^{aa}$ is skew-Hermitian and
\[\forall z\geq 0, \quad \widehat{a}^\e(\omega,z)=\textbf{T}^{\xi,\e}(\omega,z)(\widehat{a}^\e_0(\omega)).\]

\section{Time Reversal in a Waveguide}\label{sect4P2p}

\begin{figure} \begin{center}
\begin{tabular}{cc}
\includegraphics*[scale=0.27]{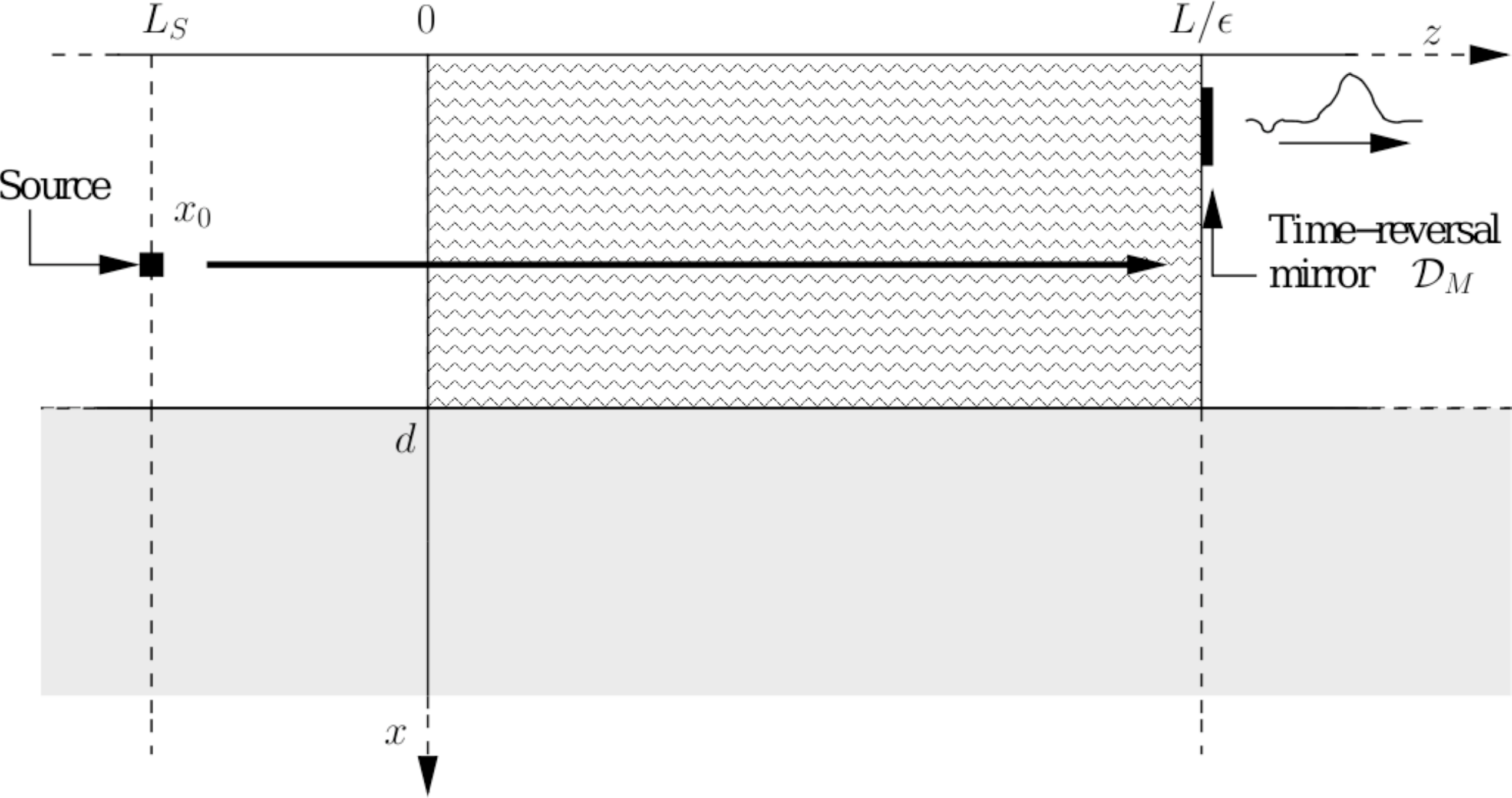} &\includegraphics*[scale=0.27]{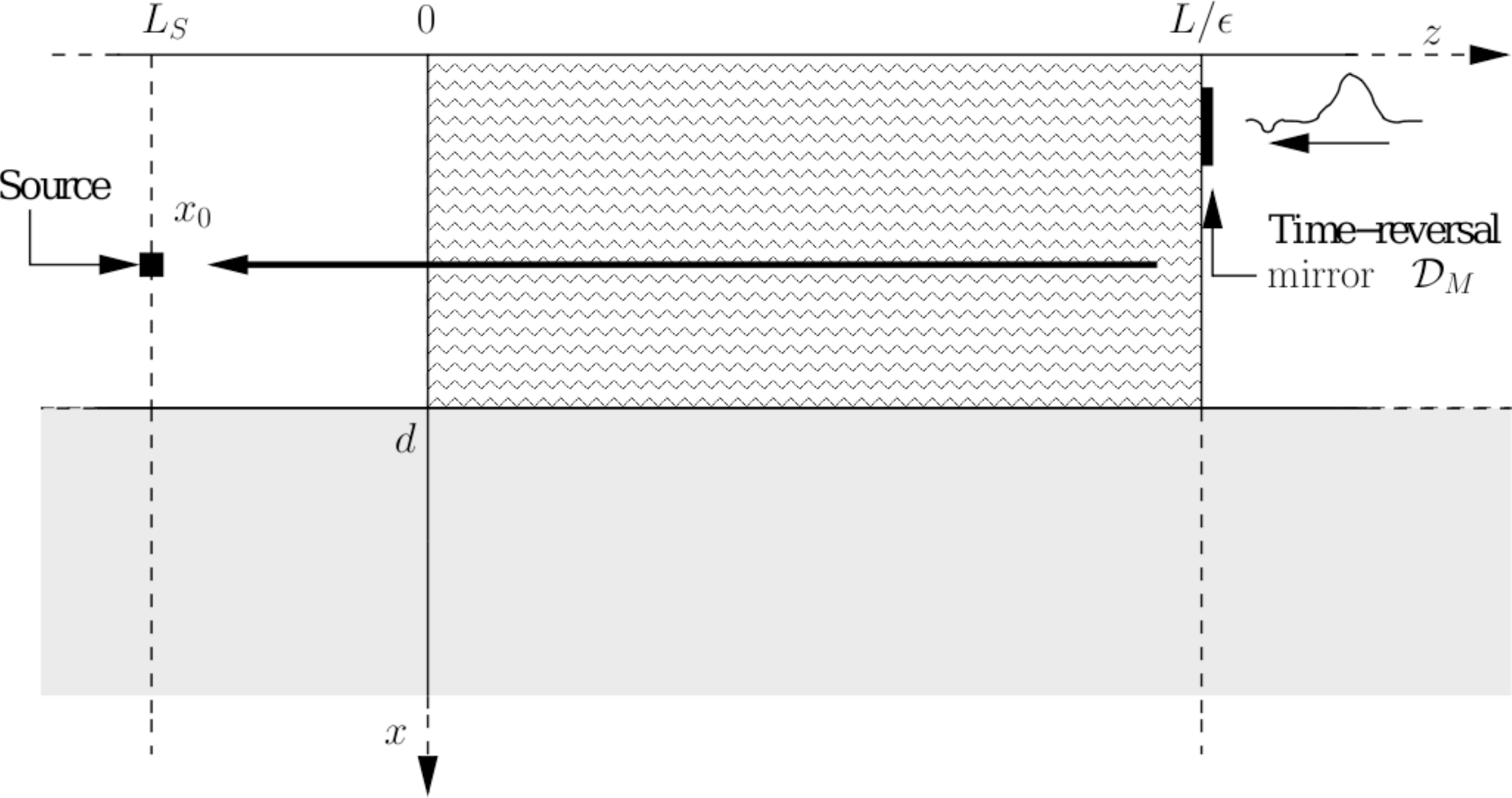}\\
$(a)$ & $(b)$
\end{tabular}
\end{center}
\caption{\label{figtrmP2p}
Representation of the time-reversal experiment. In $(a)$ we represent the first step of the experiment, and in $(b)$ we represent the second step of the experiment. 
}\end{figure}

Time-reversal experiments with sonar in shallow water \cite{kuperman,kuperman2} were carried out by William Kuperman and his group in San Diego.
This experiment is carried out in two steps. In the first step (see Figure \ref{figtrmP2p} $(a)$), a source sends a pulse into the medium. The wave propagates and is recorded by a device called a time-reversal mirror. A time-reversal mirror is a device that can receive a signal, record it, and resend it time-reversed into the medium. In other words, what is recorded first is send out last. In the second step (see Figure \ref{figtrmP2p} $(b)$), the wave emitted by the time-reversal mirror has the property of refocusing near the original source location, and it has been observed that random inhomogeneities enhance refocusing \cite{papanicolaou3,clouet,papanicolaou2,fink1,fink2,trfpsource,book,papanicolaou4}. This experiment has already been analyzed in waveguides with bounded cross-section in \cite[Chapter 20]{book} and \cite{papa,gomez}. However, in contrast with all these results where the random medium improve the refocusing, we show in this section in a context  of underwater acoustics that the random inhomogeneities deteriorate the refocusing property. In \cite{ammari} the authors observe such deterioration but in their context it is induced by viscosity effect in the wave propagation model. We show in this section that this effect is simply induced by the inhomogeneities of the propagation medium through the coupling mechanism between the propagating and the radiating modes.

\subsection{First Step of the Experiment}

\begin{figure}
\begin{center}\includegraphics*[scale=0.37]{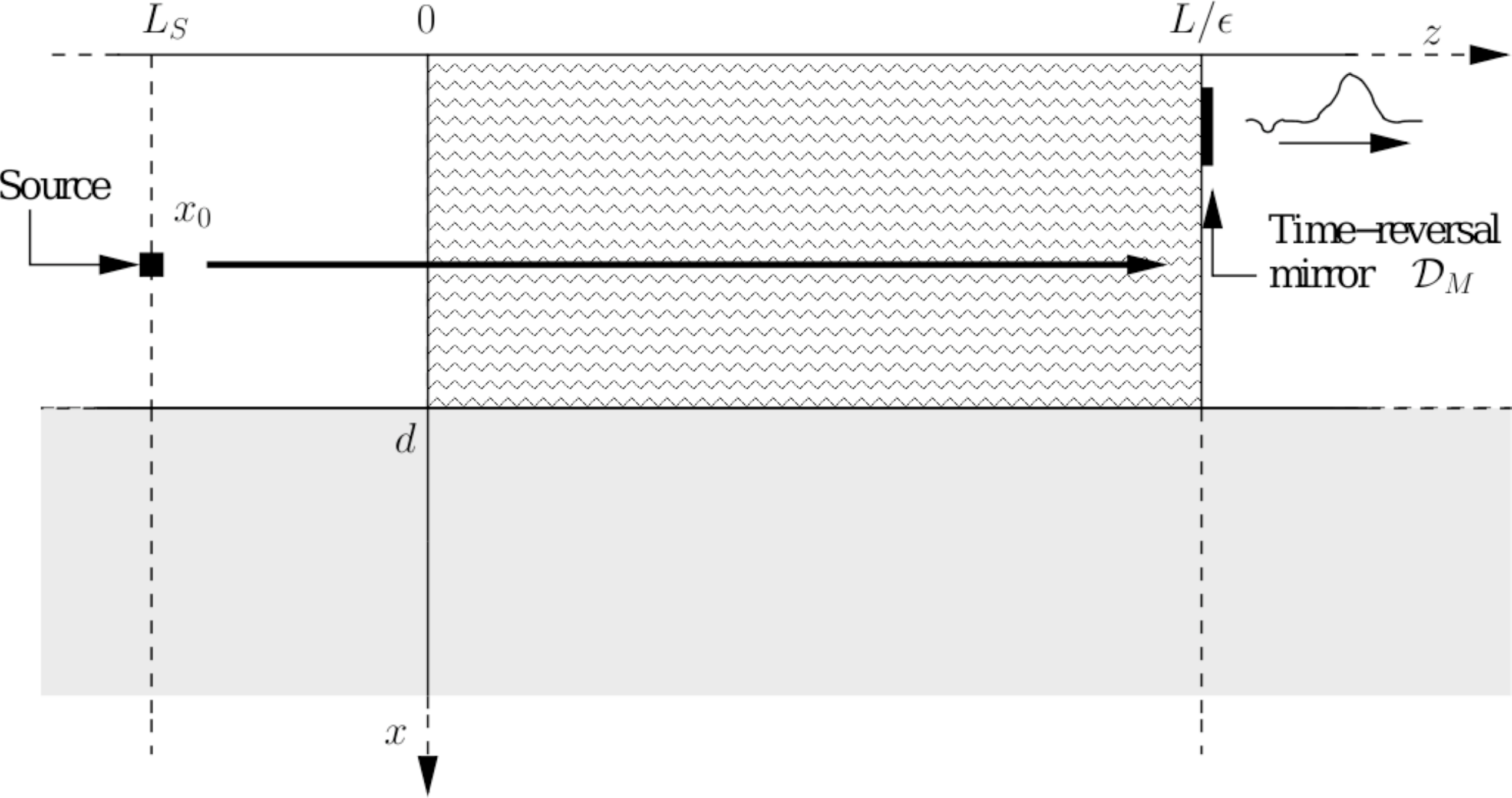}\end{center}
\caption{\label{figtrm2P2p} Representation of the first step of the time-reversal experiment.}
\end{figure}

In the first step of the experiment (see Figure \ref{figtrm2P2p}), a source sends a pulse into the medium, the wave propagates and is recorded by the time-reversal mirror located in the plane $z=L / \e$. We assume that the time-reversal mirror occupies the transverse subdomain $\D_{M}\subset [0,d]$ and in the first step of the experiment the time-reversal mirror plays the role of a receiving array. The transmitted wave is recorded for a time interval $\big[ \frac{t_0}{\e}, \frac{t_1}{\e}\big] $ and is re-emitted time-reversed into the waveguide toward the source. We have chosen such a time window because it is of the order of the total travel time of the section $[0,L/\e]$.

In this paper, the source profile $\Psi^\e ( t,x)$ is given, in the frequency domain, by
\begin{equation} \label{profsourceP2p}\begin{split}
\widehat{\Psi}^\e _q( \omega,x)=\frac{1}{\e^q}&\widehat{f}\left(\frac{\omega-\omega_0}{\e^q}\right)\\
&\times\left[\sum_{j=1}^{\N{}}\phi_j(\omega,x_0)\phi_j(\omega,x)+\int_{(-S,-\xi)\cup(\xi,\ko{})}\phi_\ga(\omega,x_0)\phi_\ga(\omega,x)d \ga\right],\end{split}\end{equation}
with $q>0$. The restriction $q>0$ allows us to freeze the number of propagating and radiating modes, introduced below, and gives simpler expressions of the transmitted field. Let us note that $S$ can be arbitrarily large and $\xi$ can be arbitrarily small, so that the transverse profile \eqref{profsourceP2p} is an approximation of a Dirac distribution at $x_0$, which models a point source at $x_0$. Moreover, $\frac{1}{\e^q}\widehat{f}(\frac{\omega-\omega_0}{\e^q})$ is the Fourier transform of $f(\e^q t)e^{-i \omega_0 t}$, which is a pulse with bandwidth of order $\e^q$ and carrier frequency $\omega_0$. In this paper, we are interested by a source emitting a broadband pulse, that is for $q\in (0,1)$. For our broadband source term with pulse width of order $1/\e^q$, smaller than the propagation distance, the propagating modes are separated in time by the modal dispersion. With this kind of source, waiting long enough to record all the train of pulses at the time-reversal mirror, one can observe a self-averaging effect on the refocused pulse. This statistical stability implies that the refocused pulse does not depend on the particular realization of the random medium. This phenomenon has been widely studied in different contexts and there are many references about it \cite{papanicolaou3,bal2,papanicolaou2,fink1,book,papa}. The case $q=1$, that we do not treat in this paper corresponds to the narrowband case. In this case the order of the pulse width is comparable to the propagation distance. Consequently, the modes overlap during the propagation and then the statistical stability of the time-reversal experiment depends on the number of propagating modes \cite{book,papa}. However, for the sake of simplicity, we consider only the case $q=1/2$ but the following analysis can be carried out $\forall q\in(0,1)$. 

According to \cite{book,papa,gomez2}, the evanescent part of the wave field decreases exponentially fast with the propagation distance. For more convenient manipulations in the study of the time-reversal experiment we assume that the source location $L_S$ is sufficiently far away from $0$ so that the evanescent modes generated by the source are negligible. With this assumption and using \eqref{initfieldP2}, we can assume that the incident pulse coming from the left is given, at $z=0$, by: 
\[p^{\xi,\e}_{inc}(t,x,0)=\frac{1}{2\pi}\int \left[\sum_{j=1}^{\N{}} \frac{ \widehat{a}^\e_{j,0} (\omega) }{\sqrt{\Bh{j}{}}} \phi_j(\omega,x)+\int_{\xi}^{\ko{}} \frac{\widehat{a}^\e_{\ga,0} (\omega)}{\ga^{1/4}}\phi_\ga(\omega,x)d \ga\right]e^{-i\omega t}d\omega,\]    
where
\begin{equation}\label{initcondP21p}
\widehat{a}^\e_{j,0} (\omega)  = \frac{\sqrt{\Bh{j}{}}}{2\e^q} \widehat{f}\left(\frac{\omega-\omega_0}{\e^q}\right)\phi_j(\omega,x_0)e^{-i\Bh{j}{}L_S}=\frac{1}{2\e^q} \widehat{f}\left(\frac{\omega-\omega_0}{\e^q}\right)\tilde{a}_j(\omega)\end{equation}
$\forall j \in \big\{1,\dots,\N{} \big\}$,
\begin{equation}\label{initcondP22p}
\widehat{a}^\e_{\ga,0} (\omega)=\frac{\ga ^{1/4}}{2\e^q} \widehat{f}\left(\frac{\omega-\omega_0}{\e^q}\right)\phi_\ga(\omega,x_0) e^{-i \sga L_S} = \frac{1}{2\e^q} \widehat{f}\left(\frac{\omega-\omega_0}{\e^q}\right)\tilde{a}_\ga(\omega)\end{equation}
for almost every $\ga \in(\xi, \ko{})$. Let us remark that this assumption is not restrictive and all the results of this paper are valid for any $L_S<0$. Indeed, according to Proposition 4.2 in \cite{gomez2}, in the asymptotic $\e\to 0$, the information about the evanescent part of the source profile are lost during the propagation in the random section $[0,L/\e]$, and therefore they play no role in the pulse propagation and in the time-reversal experiment. An efficient way to do not loss the information about the evanescent part of the source term has been developed in \cite{science} and studied in \cite{gomez}.

Finally, according to Section \ref{sect2P2p}, the wave recorded by the time-reversal mirror is given by
\[\begin{split}
p_{tr}\left(t,x,\frac{L}{\e}\right)=\frac{1}{4\pi \se}&\int \widehat{f}\left(\frac{\omega-\omega_0}{\se}\right)\\
&\quad \times \left[ \sum_{j=1}^{\N{}} \frac{1}{\sqrt{\Bh{j}{}}}\textbf{T}_j^{1,\xi,\e}(\omega,L)(\tilde{a} (\omega)) \phi_j(\omega,x)e^{i\Bh{j}{}\frac{L}{\e}}e^{-i\omega t}\right.\\
&\quad + \left.\int_{\xi}^{\ko{}} \frac{1}{\ga^{1/4}}\textbf{T}_\ga^{1,\xi,\e}(\omega,L)(\tilde{a} (\omega))\phi_\ga(\omega,x)e^{i\sga\frac{L}{\e}}d \ga e^{-i\omega t} \right]d\omega,
\end{split}\]
where $\textbf{T}^{\xi,\e}(\omega,L)$ is the transfer operator solution of \eqref{transferP2p}.

\subsection{Second Step of the Experiment}

\begin{figure}
\begin{center}\includegraphics*[scale=0.37]{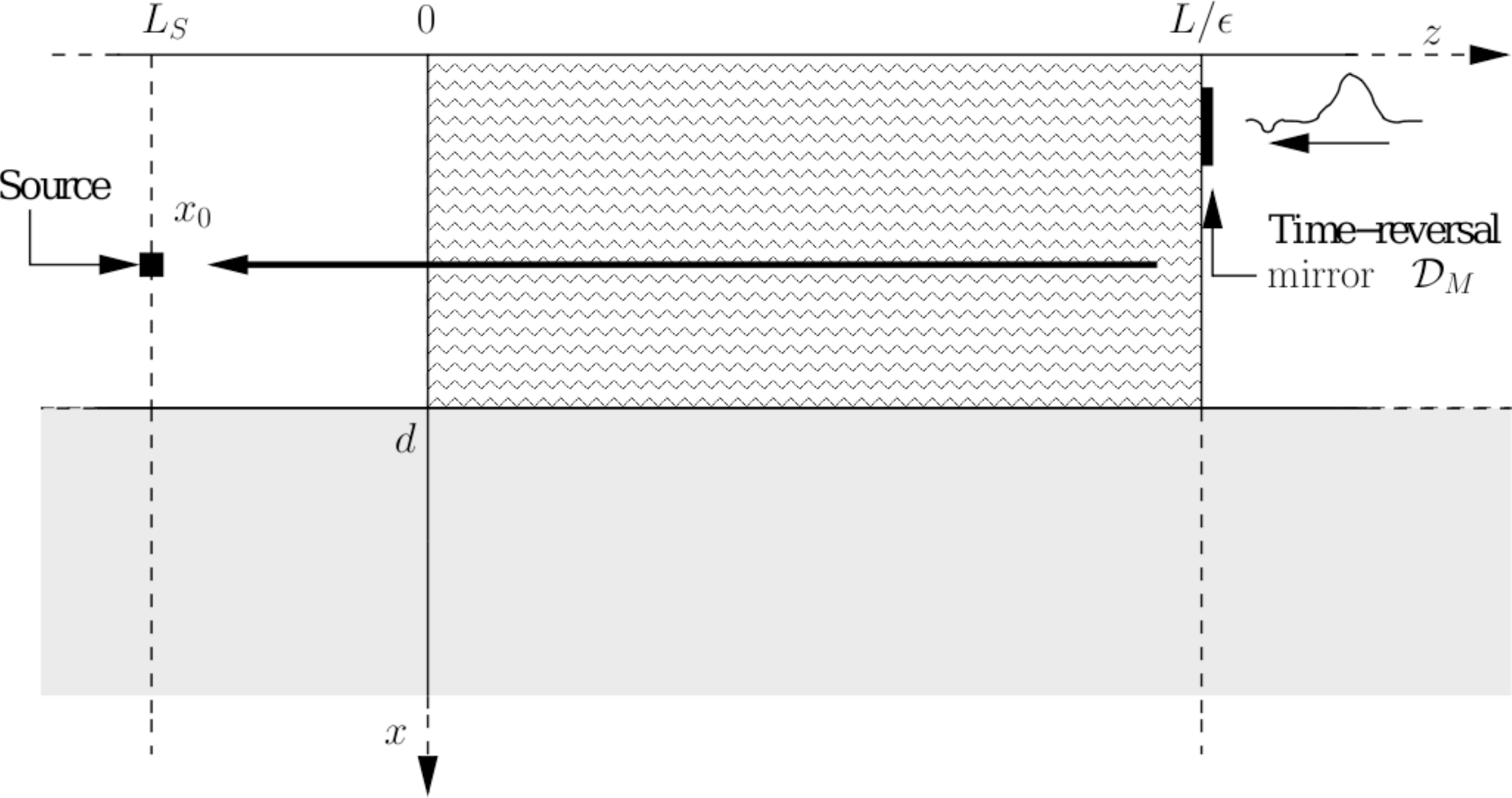}\end{center}
\caption{\label{figtrm3P2p} Representation of the second step of the time-reversal experiment.}
\end{figure}

In the second step of the experiment (see Figure \ref{figtrm3P2p}), the time-reversal mirror plays the role of a source array, and the time-reversed signal is transmitted back. Now, the source term is given by
\[\textbf{F}^\e _{TR}(t,x,z)=-f^\e _{TR}(t,x)\delta(z-L / \e)\textbf{e}_{z},\] 
with
\[f^\e _{TR}(t,x)=\Mt{\frac{t_1}{\e}- t}{\frac{L_M}{\e}}G_1 (t_1-\e t) G_2 (x),\]
where
\[G_1 (t)=\1 _{[t_0 ,t_1]} (t) \quad \text{and} \quad G_2 (x)=\1_{\D_M}(x).\]
Here, $G_1$ represents the time window in which the transmitted wave is recorded, and $G_2$  represents the spatial window in which the transmitted wave is recorded. In our study, we are interested in the spatial effects of the refocusing, so we assume that we record the field for all time at the time-reversal mirror, that is the source has the form 
\[
f^\e _{TR}(t,x)=p_{tr}\left(\frac{t_1}{\e}-t,x,\frac{L}{\e}\right)G_2 (x).
\]
Now, we are interested in the propagation from $z=L /\epsilon$ to $z=0$. The decomposition with respect to the resolution of the identity $\Pi_\omega$ associated to $R(\omega)$ (see Section \ref{spectralP2}) gives
\[\widehat{p}_{TR}(\omega,x,z)=\sum_{m=1}^{\N{}} \frac{\widehat{b}_{m}(\omega,z)}{\sqrt{\Bh{m}{}}} e^{-i\Bh{m}{}z}\phi_m(\omega,x)+\int_{\xi}^{k^2(\omega)} \frac{\widehat{b}_{\gamma}(\omega,z)}{\gamma^{1/4}} e^{-i\sga z}\phi_\ga(\omega,x)d \ga,\]
with
\begin{equation*}\begin{split}
\widehat{b}_{m}(\omega,L)&=\frac{\sqrt{\Bh{m}{}}}{2}e^{i\Bh{m}{}\frac{L}{\e}}\big< \widehat{f}^\e _{TR}(\omega,.),\phi_m(\omega,.) \big>_{H}, \\
\widehat{b}_{\ga}(\omega,L)&=\frac{\ga^{1/4}}{2}e^{i\sga\frac{L}{\e}}\big< \widehat{f}^\e _{TR}(\omega,.),\phi_\ga(\omega,.) \big>_{H}
\end{split}\end{equation*}
in $\esp_\xi$. Then, at the source location $z=L_S$, we obtain
\[\widehat{p}_{TR}(\omega,x,L_S)=\sum_{n=1}^{\N{}}\frac{\widehat{b}_n(\omega,0)}{\sqrt{\Bh{n}{}}}e^{i\Bh{n}{}L_S}\phi_n(\omega,x)+\int_{\xi}^{k^2(\omega)}\frac{\widehat{b}_\ga (\omega,0)}{\ga^{1/4}}e^{i\sga L_S}\phi_\ga(\omega,x)d\ga,\]
with
\begin{equation*}
\widehat{b}(\omega,0)=\overline{\big(\textbf{T}^{\xi,\e}\big)^{\ast}(\omega,L)}\big(\widehat{b}(\omega,L)\big),
\end{equation*}
where $\big(\textbf{T}^{\xi,\e}\big)^{\ast}(\omega,z)$ stands for the adjoint operator of $\textbf{T}^{\xi,\e}(\omega,z)$. 
Consequently,  one can write
\[\widehat{p}_{TR}(\omega,x,L_S)=\big<\textbf{T}^{\xi,\e}(\omega,L)(\tilde{b}_x(\omega)),\overline{\widehat{b}(\omega,L)}\big>_{\esp_\xi },\]
where
\begin{equation}\label{bcoeftrP2p}\tilde{b}_{x,n}(\omega)=\frac{1}{\sqrt{\Bh{n}{}}}\phi_n(\omega,x) e^{-i\Bh{n}{}L_S}\quad \text{and}\quad \tilde{b}_{x,\ga}(\omega)=\frac{1}{\ga^{1/4}}\phi_\ga(\omega,x) e^{-i\sga L_S},\end{equation}
and
\[\begin{split}
\widehat{b}^2_m(\omega,L)&=\frac{1}{4\se}\overline{\widehat{f}\left(\frac{\omega-\omega_0}{\se}\right)}e^{i\omega t_1}\big<\overline{\textbf{T}^{1,\xi,\e}(\omega,L)(\tilde{a}(\omega))},\lambda^\e_m(\omega)\big>_{\esp_\xi},\\
\widehat{b}^2_\ga (\omega,L)&=\frac{1}{4\se}\overline{\widehat{f}\left(\frac{\omega-\omega_0}{\se}\right)}e^{i\omega t_1}\big<\overline{\textbf{T}^{1,\xi,\e}(\omega,L)(\tilde{a}(\omega))},\lambda^\e_\ga (\omega)\big>_{\esp_\xi},
\end{split}\]
in $\esp_\xi$, where $\lambda^\e(\omega)$ is defined by
\begin{equation}\label{lambdacoeftrP2p}\begin{split}
\lambda^\e(\omega)_{mj}&=\sqrt{\frac{\Bh{m}{}}{\Bh{j}{}}}e^{-i(\Bh{m}{}-\Bh{j}{})\frac{L}{\e}}M_{mj}(\omega),\\
\lambda^\e(\omega)_{m\ga'}&=\sqrt{\frac{\Bh{m}{}}{\sgap}}e^{-i(\Bh{m}{}-\sgap)\frac{L}{\e}}M_{m\ga'}(\omega),\\
\lambda^\e(\omega)_{\ga j}&=\sqrt{\frac{\sga}{\Bh{j}{}}}e^{-i(\sga-\Bh{j}{})\frac{L}{\e}}M_{\ga j}(\omega),\\
\lambda^\e(\omega)_{\ga \ga'}&=\frac{\ga^{1/4}}{{\ga'}^{1/4}}e^{-i(\sga-\sgap)\frac{L}{\e}}M_{\ga \ga'}(\omega),\\
\end{split}\end{equation}
and with
\[M_{rs}(\omega)=\int_0^d G_2(x)\phi_r(\omega,x)\phi_s(\omega,x)dx \]
for $(r,s)\in\big(\{1,\dots,\N{}\}\cup(\xi,\ko{})\big)^2$. $(M_{rs}(\omega))$ represents the coupling produced by the time-reversal mirror between the modes during the two steps of the time-reversal experiment. Consequently,
\[
\widehat{p}_{TR}(\omega,x,L_S)=\frac{1}{4\se}\overline{\widehat{f}\left(\frac{\omega-\omega_0}{\se}\right)}e^{i\omega t_1}\big<\textbf{U}^{\xi,\e}(\omega,L)\big(\tilde{a}(\omega),\tilde{b}_x(\omega)\big),\lambda^\e(\omega)\big>_{\esp_\xi \otimes \esp_\xi}.
\]
Here, we consider the tensorial space $\esp_\xi \otimes \esp_\xi=\left\{ \lambda \otimes \mu,\quad (\lambda,\mu)\in(\esp_\xi)^2 \right\}$, with
$(\lambda \otimes \mu )_{rs}=\lambda_r \mu_s$ for $(r,s)\in\big(\{1,\dots,\N{}\}\cup(\xi,\ko{})\big)^2$ and $\forall (\lambda,\mu)\in(\esp_\xi)^2$.
This space is equipped with the inner product defined by
\[\begin{split}
\big<\lambda,\mu\big>_{\esp_\xi \otimes \esp_\xi}&=\sum_{j,l=1}^{\N{}}\lambda_{jl}\overline{\mu_{jl}}+\sum_{j=1}^{\N{}}\int_{\xi}^{\ko{}}\lambda_{j\ga'}\overline{\mu_{j\ga'}}d\ga'\\
&\quad+\int_{\xi}^{\ko{}}\sum_{l=1}^{\N{}}\lambda_{\ga l}\overline{\mu_{\ga l}}d\ga+\int_{\xi}^{\ko{}}\int_{\xi}^{\ko{}}\lambda_{\ga \ga'}\overline{\mu_{\ga \ga'}}d\ga d\ga'
\end{split}\]
 $\forall (\lambda,\mu)\in(\esp_\xi \otimes \esp_\xi)^2$. Finally, the time-reversal kernel $\textbf{U}^{\xi,\e}(\omega,L)$ is defined by 
 \begin{equation}\label{TRkernel}\textbf{U}^{\xi,\e}(\omega,L)(y^1,y^2)=\overline{\textbf{T}^{\xi,\e}(\omega,L)(y^1) }\otimes \textbf{T}^{\xi,\e}(\omega,L)(y^2)\end{equation}
$\forall (y^1,y^2)\in(\esp_\xi)^2$, describing the two steps of the time-reversal experiment through the random medium thanks to the transfer operator  $\textbf{T}^{\xi,\e}(\omega,L)$ satisfying \eqref{transferP2p}. 

We study the refocused wave in a time window of order $1/\se$ comparable to the pulse width, and centered at time $t_{obs}/\e$, which is of the order the total travel time for a distance of order $1/\e$. 
Consequently, we will study the refocusing of the refocused wave at the source location $z=L_S$ given by
\begin{equation}\label{refocwaveform}\begin{split}
p_{TR}\Big(&\frac{t_{obs}}{\e}+\frac{t}{\se},x,L_S\Big)=\frac{1}{2\pi}\int\widehat{p}_{TR}(\omega,x,L_S)e^{-i\omega t}d\omega\\
&=\frac{1}{8\pi\se}\int \overline{\widehat{f}\left(\frac{\omega-\omega_0}{\se}\right)}\big<\textbf{U}^{\xi,\e}(\omega,L)\big(\tilde{a}(\omega),\tilde{b}_x(\omega)\big),\lambda^\e(\omega)\big>_{\esp_\xi \otimes \esp_\xi}e^{i\omega \Big(\frac{t_1-t_{obs}}{\e}-\frac{t}{\se}\Big)}d\omega,
\end{split}\end{equation}
where $\tilde{a}(\omega)$ is defined by \eqref{initcondP21p} and \eqref{initcondP22p}, $\tilde{b}_x(\omega)$ is defined by \eqref{bcoeftrP2p}, and $\lambda^\e(\omega)$ is defined by \eqref{lambdacoeftrP2p}.

In what follows, we consider a time-reversal mirror of the form $\D_M=[d_1,d_2]$ with 
\[d_2 = d_M +\lambda_{oc} ^{\alpha_M} \td{2} \text{ and } d_1=d_M -\lambda_{oc} ^{\alpha_M}\td{1},\]
where $d_M \in (0,d)$, $(\td{2},\td{1})\in (0,+\infty)^2$, and $\alpha_M \in [0,1]$. 
Here, $\lambda_{oc} =2\pi c/(n_1\omega_0)$ is the carrier wavelength in the ocean section $[0,d]$ of the waveguide. The time-reversal coupling matrix are therefore given by  
 \begin{equation}\label{matrixmirror}\begin{split} M_{jl}(\omega)&=(d_2 -d_1)A_j(\omega)A_l(\omega)\\
 &\quad\times\left[ \cos\left((\sigma_j(\omega)-\sigma_l(\omega))\frac{d_2 +d_1}{2d}\right)\textrm{sinc} \left((\sigma_j(\omega)-\sigma_l(\omega))\frac{d_2 -d_1}{2d}\right) \right. \\ 
 &\left. \quad \quad-\cos\left((\sigma_j(\omega)+\sigma_l(\omega))\frac{d_2 +d_1}{2d}\right)\textrm{sinc} \left((\sigma_j(\omega)+\sigma_l(\omega))\frac{d_2 -d_1}{2d}\right) \right],\end{split}\end{equation}
for $(j,l)\in\{1,\dots,\N{}\}^2$, where $A_j(\omega)$ and $\sigma_j(\omega)$ are defined in Section \ref{spectralP2}. We give only the coefficients $M_{jl}(\omega)$ for $(j,l)\in\{1,\dots,\N{}\}^2$, because in what follows only these terms will play a role. 
The parameter $\alpha_M$ represents the order of the magnitude of the size of the time-reversal mirror with respect to the wavelength in the ocean cross-section $[0,d]$.  In fact, we will see that the size of the mirror plays a role in the homogeneous case only when it is of the order the carrier wavelength $\lambda_{oc}=2\pi c/(n_1 \omega_0)$. 

In the following section we study the transverse profile of the refocused wave in the continuum limit $\N{_0} \gg 1$ of a large number of propagating modes, which corresponds to the regime $\omega_0 \nearrow +\infty$. However, we know that the main focal spot must be of order $\lambda_{oc}$, which tends to $0$ in this continuum limit $\N{_0}\gg1$. Consequently, we study the transverse profile of the refocused wave in a spatial window of size $\lambda_{oc}$ centered around $x_0$.

\subsection{Refocused Field in a Homogeneous Waveguide}\label{rfhwintP2p}

To understand what are the effects produced by the random perturbations of the propagation medium on the time-reversal experiment, we study first the refocused wave obtained in a homogeneous waveguide.
Let us consider the refocused wave in a time window of order $1/\se$, which is comparable to the pulse width, and centered at time $t_{obs}/\e$, which is of the order the total travel time for a distance of order $1/\e$.
In this section we assume that the medium is homogeneous, so that $\textbf{T}^{\xi,\e}(\omega,L)=Id$. Then, the refocused wave at the original source location is given by
\[ \begin{split} 
p_{TR}\Big(\frac{t_{obs}}{\e}+\frac{t}{\se},x,&L_S\Big)=e^{i\omega_0\frac{t_1 -t_{obs}}{\e}}e^{-i\omega_0\frac{t}{\se}}\cdot\frac{1}{4}\sum_{j,m=1}^{\N{_0}}e^{i(\Bh{m}{_0}-\Bh{j}{_0})\left(-L_S+\frac{L}{\e}\right)}M_{jm}(\omega_0)\\
&\times \phi_j(\omega_0,x_0)\phi_m(\omega_0,x) K^{\omega_0}_{j,m,L}\ast f\left(\frac{(\beta'_m(\omega_0)-\beta'_j(\omega_0))L+t_1-t_{obs}}{\se}-t\right)\\
&+\mathcal{O}(\se),
\end{split}\]
where
\begin{equation}\label{kernelP2p}
\widehat{K^{\omega_0}_{j,m,L}}(\omega)=\widehat{K^{\omega_0}_{j,L}}(\omega)\overline{\widehat{K^{\omega_0}_{m,L}}(\omega)}=e^{i(\beta''_j(\omega_0)-\beta''_m(\omega_0))L\frac{\omega^2}{2}},
\end{equation}
and $K^{\omega_0}_{j,j,L}=\delta_0$. Consequently, in the asymptotic $\e\to 0$, we can observe a refocused wave only for a finite set of times given by
\begin{equation}\label{timetrP2p}
t_{jm}=t_1+(\beta'_m(\omega_0)-\beta'_j(\omega_0))L.
\end{equation}
For $m\not=j$, we obtain
\[ \begin{split} 
p_{TR}\Big(\frac{t_{jm}}{\e}+\frac{t}{\se},x,L_S\Big)=&e^{i\omega_0\frac{t_1 -t_{jm}}{\e}}e^{-i\omega_0\frac{t}{\se}}e^{i(\Bh{m}{_0}-\Bh{j}{_0})\left(-L_S+\frac{L}{\e}\right)}M_{jm}(\omega_0)\\
&\times \frac{1}{4}\phi_j(\omega_0,x_0)\phi_m(\omega_0,x) K^{\omega_0}_{j,m,L}\ast f(-t)\\
&+\mathcal{O}(\se).
\end{split}\]
At time $t_{jm}$ ($j\not=m$) one can observe only the $m$th mode, emitted by the time-reversal mirror during the second step of the experiment, coupled with the $j$th modes recorded by the time-reversal mirror during the first step. This coupling is produced by the time-reversal mechanism through the time-reversal mirror and  characterized by the coupling matrix $M_{jm}(\omega_0)$. Moreover, let us note that the refocused wave shape is dispersed by the kernel $K^{\omega_0}_{j,L}(t)$ during the first step of the experiment and by $K^{\omega_0}_{m,L}(-t)$ during the second step.

Now, for $t_{obs}=t_1$ we obtain
\[ p_{TR}\Big(\frac{t_{1}}{\e}+\frac{t}{\se},x,L_S\Big)=e^{-i\omega_0\frac{t}{\se}}f(-t)H^{\alpha_M}_{x_0}(\omega_0,x)+\mathcal{O}(\se),\]
where
\[H^{\alpha_M}_{x_0}(\omega_0,x)=\frac{1}{4}\sum_{j=1}^{\N{_0}}M_{jj}(\omega_0) \phi_j(\omega_0,x_0)\phi_j(\omega_0,x).\]
Here, we have a contribution of all the propagating modes. The refocused wave is a superposition of modes where each mode is coupled with itself by the time-reversal mirror through the terms $M_{jj}(\omega_0)$. We describe in the following proposition the transverse profile of the time-reversed pulse in a very simple way in the regime $\omega_0\nearrow +\infty$, which correspond to the continuum limit $\N{_0}\gg1$ of a large number of propagating modes.

\begin{prop}\label{prop17P2p} For $\alpha_M\in[0,1)$, the transverse profile of the refocused wave in the continuum limit  $\N{_0}\gg1$ is given by
\[
\lim_{\omega_0 \to +\infty}\frac{\lambda_{oc}^{1-\alpha_M}}{\theta}H^{\alpha_M}_{x_0}\Big(\omega_0,x_0+\frac{\lambda_{oc}}{\theta}\tilde{x}\Big)=\frac{\tilde{d}_2+\tilde{d}_1}{d}\emph{sinc}(2\pi\tilde{x}).
\]
The width of the focal spot is therefore given by the diffraction limit $\lambda_{oc}/(2\theta)$.
\end{prop}
The proof of Proposition \ref{prop17P2p} is given in Section \ref{proof0} 	
As a result the sinc function describes the asymptotic transverse profile of the refocused wave in the continuum limit.

The sinc profile has already been obtained in different contexts in time reversal to describe the refocused transverse profile. Moreover, the size of the focal spot is of order $\lambda_{oc}/(2\theta)$, where $\lambda_{oc}$ is the carrier wavelength of the ocean section of the waveguide, and $\theta$ depends on the contrast of the refractive index between the ocean section and the bottom of the waveguide. In \cite{papanicolaou2,book,papa,gomez} for instance the random perturbations of the medium improve the time-reversal refocusing, but we will see in what follows that this statement is no more true in our context if the mode coupling mechanism between the propagating and the radiating mode is not negligible.

\subsection{Limit Theorem}

To describe the effects of the random medium on the time-reversed wave \eqref{refocwaveform} we need to know the asymptotic distribution of the process  $\textbf{U}^{\xi,\e}(\omega,.)$, defined by \eqref{TRkernel} as $\e$ goes to $0$ and $\xi$ goes to $0$. First, let us remark that $\forall (y^1,y^2)\in(\esp_\xi)^2$, with $\esp_\xi=\mathbb{C}^{\N{}}\times L^2(\xi,\ko{})$,
\[\|\textbf{U}^{\xi,\e} (\omega,z)(y^1,y^2)\|^2_{\esp_{\xi}\otimes \esp_\xi}=\|y^1\otimes y^2\|^2_{\esp_{\xi}\otimes \esp_\xi} \quad \forall z\geq 0, \]
and then let us introduced some notations.
Let $r_y=\|y^1\otimes y^2\|_{\esp_{\xi}\otimes \esp_\xi}$,
\[\mathcal{B}_{r_y, \esp_{\xi}\otimes \esp_\xi}=\left\{\lambda \in \esp_\xi\otimes \esp_\xi, \|\lambda\|_{\esp_\xi\otimes \esp_\xi } \leq r_y\right\}\] 
the closed ball with radius $r_y$, and $\{g_n, n\geq 1\}$ a dense subset of $\mathcal{B}_{r_y,\esp_\xi\otimes \esp_\xi}$. We equip  $\mathcal{B}_{r_y, \esp_\xi\otimes \esp_\xi}$ with the distance $d_{\mathcal{B}_{r_y, \esp_\xi\otimes \esp_\xi}}$ defined by
\[d_{\mathcal{B}_{r_y,\esp_\xi\otimes \esp_\xi}}(\lambda, \mu)=\sum_{j=1}^{+\infty}\frac{1}{2^j}\left\lvert\big<\lambda-\mu,g_n\big>_{\esp_\xi\otimes \esp_\xi}\right\rvert\]
$\forall (\lambda,\mu)\in{(\mathcal{B}_{r_y,\esp_\xi\otimes \esp_\xi})}^2$, so that $(\mathcal{B}_{\esp_\xi} ,d_{\mathcal{B}_{r_y, \esp_\xi\otimes \esp_\xi}})$ is a compact metric space.

In the following theorem, we give only the drifts of the infinitesimal generators because only this part is of interest in what follows. 

\begin{thm}\label{thasympP2p1}
$\forall (y^1,y^2)\in( \esp_\xi)^2$, the stochastic process $\emph{\textbf{U}}^{\xi,\e} (\omega,.)(y^1,y^2)$ converges in distribution on $\mathcal{C}([0,+\infty),(\mathcal{B}_{r_y,\esp_\xi \otimes \esp_\xi},d_{\mathcal{B}_{r_y,\esp_\xi \otimes \esp_\xi}}))$ as $\e\to 0$ to a limit denoted by $\emph{\textbf{U}}^{\xi} (\omega,.)(y^1,y^2)$, unique solution of a well-posed martingale problem on $\esp_\xi \otimes \esp_\xi$ starting from $y^1\otimes y^2$. Moreover, $\forall (y^1,y^2)\in( \esp_0)^2$, the stochastic process $\emph{\textbf{U}}^{\xi} (\omega,.)(y^1,y^2)$ converge in distribution on $\mathcal{C}([0,+\infty),(\mathcal{B}_{r_y,\esp_0 \otimes \esp_0},d_{\mathcal{B}_{r_y,\esp_0 \otimes \esp_0}}))$ as $\xi \to 0$ to a limit denoted by $\emph{\textbf{U}}^0 (\omega,.)(y^1,y^2)$. This limit is the unique solution of the well-posed martingale problem on $\esp_0 \otimes \esp_0$ starting from $y^1\otimes y^2$, and with drift given by  
\[\mathcal{L}^\omega_{1}+\mathcal{L}^\omega_{2},\]
where
\[\begin{split}
\mathcal{L}^\omega_{1}=&\sum_{\substack{j,l=1\\j\not=l}}^{\N{}}\tilde{\Gamma}^c_{jl}(\omega)\big(U_{ll}\partial_{U_{jj}}+\overline{U_{ll}}\partial_{\overline{U_{jj}}}\big)\\
&+\frac{1}{2}\sum_{j,l=1}^{\N{}}\big[\Gamma^c_{jj}(\omega)+\Gamma^c_{ll}(\omega)-\big(\Gamma^1_{jj}(\omega)+\Gamma^1_{ll}(\omega)-2\tilde{\Gamma}^1_{jl}(\omega)\big)\big]\big(U_{jl}\partial_{U_{jl}}+\overline{U_{jl}}\partial_{\overline{U_{jl}}}\big)\\
&+\frac{1}{2}\sum_{j=1}^{\N{}}\int_{\xi}^{\ko{}}\big[\Gamma^c_{jj}(\omega)-\Gamma^1_{jj}(\omega)\big]\big(U_{j\ga_2}\partial_{U_{j\ga_2}}+\overline{U_{j\ga_2}}\partial_{\overline{U_{j\ga_2}}}\big)d\ga_2\\
&+\frac{1}{2}\int_{\xi}^{\ko{}}\sum_{l=1}^{\N{}}\big[\Gamma^c_{ll}(\omega)-\Gamma^1_{ll}(\omega)\big]\big(U_{\ga_1l}\partial_{U_{\ga_1 l}}+\overline{U_{\ga_1 l}}\partial_{\overline{U_{\ga_1 l}}}\big)d\ga_1\\
&+ \frac{i}{2} \sum_{j,l=1}^{\N{}}\big[\Gamma^s_{ll}(\omega)-\Gamma^s_{jj}(\omega)\big]\big(U_{jl}\partial_{U_{jl}}-\overline{U_{jl}}\partial_{\overline{U_{jl}}}\big)\\
&-\frac{i}{2}\sum_{j=1}^{\N{}}\int_{\xi}^{\ko{}}\Gamma^s_{jj}(\omega)\big(U_{j\ga_2}\partial_{U_{j\ga_2}}-\overline{U_{j\ga_2}}\partial_{\overline{U_{j\ga_2}}}\big)d\ga_2\\
&+\frac{i}{2}\int_{\xi}^{\ko{}}\sum_{l=1}^{\N{}}\Gamma^s_{ll}(\omega)\big(U_{\ga_1l}\partial_{U_{\ga_1 l}}-\overline{U_{\ga_1 l}}\partial_{\overline{U_{\ga_1 l}}}\big)d\ga_1,
\end{split}\]
and
\[\begin{split}
\mathcal{L}^\omega_{2,\xi}=&-\frac{1}{2}\sum_{j,l=1}^{\N{}}\big[\Lambda^{c}_j(\omega)+\Lambda^{c}_l(\omega)\big]\big(U_{jl}\partial_{U_{jl}}+\overline{U_{jl}}\partial_{\overline{U_{jl}}}\big)\\
&-\frac{i}{2}\sum_{j,l=1}^{\N{}}\big[\Lambda^{s}_l(\omega)-\Lambda^{s}_j(\omega)\big]\big(U_{jl}\partial_{U_{jl}}-\overline{U_{jl}}\partial_{\overline{U_{jl}}}\big)\\
&-\frac{1}{2}\sum_{j=1}^{\N{}}\int_\xi^{\ko{}}\big[\Lambda^{c}_j(\omega)-i\Lambda^{s}_j(\omega)\big]U_{j\ga_2}\partial_{U_{j\ga_2}}+\big[\Lambda^{c}_j(\omega)+i\Lambda^{s}_j(\omega)\big]\overline{U_{j\ga_2}}\partial_{\overline{U_{j\ga_2}}}\\
&-\frac{1}{2}\int_\xi^{\ko{}}\sum_{l=1}^{\N{}}\big[\Lambda^{c}_l(\omega)+i\Lambda^{s}_l(\omega)\big]U_{\ga_1 l}\partial_{U_{\ga_1 l}}+\big[\Lambda^{c}_l(\omega)-i\Lambda^{s}_l(\omega)\big]\overline{U_{\ga_1 l}}\partial_{\overline{U_{\ga_1 l}}}.
\end{split}\]
\end{thm}
The proof of Theorem \ref{thasympP2p1} uses exactly the same techniques as the ones developed in \cite{gomez2, these}. They are based on the perturbed-test-function method introduced in \cite{kushner} and martingale techniques. Here, we have considered the complex derivative with the following notations. If $U=U^1+iU^2\in\esp_0\otimes \esp_0$, we have $(U^1,U^2)\in(\mathcal{G}^\omega_0\otimes \mathcal{G}^\omega_0)^2$, where $\mathcal{G}^\omega_0=\mathbb{R}^{\N{}}\times L^2(0,\ko{})$. Then, the operators $\partial_U=(\partial_{U_{r,s}})$ and $\partial_{\overline{U}}=(\partial_{\overline{U_{r,s}}})$ are defined by
\[\partial_U=\frac{1}{2}(\partial_{U^1}-i\partial_{U^2})\quad\text{and} \quad \partial_{\overline{U}}=\frac{1}{2}(\partial_{U^1}+i\partial_{U^2}),\]
with $\forall f \in \mathcal{C}^1((\mathcal{G}^\omega_0\otimes \mathcal{G}^\omega_0)^2,\mathbb{R})$ and $\forall \lambda=(\lambda^1,\lambda^2)\in (\mathcal{G}^\omega_0\otimes \mathcal{G}^\omega_0)^2$
\[\begin{split}
 \sum_{n=1,2}&\Big[\sum_{j,l=1}^{\N{}}\lambda^n_{jl}\partial_{U^n_{jl}}f(v^1,v^2) +\sum_{j=1}^{\N{}} \int_\xi^{\ko{}}
 \lambda^n_{j\ga_2}\partial_{U^n_{j\ga_2}}f(v^1,v^2)d\ga_2\\
&+\int_\xi^{\ko{}} \sum_{l=1}^{\N{}}\lambda^n_{\ga_1l}\partial_{U^n_{\ga_1 l}}f(v^1,v^2) d\ga_1+\int_\xi^{\ko{}}\int_\xi ^{\ko{}}\lambda^n_{\ga_1 \ga_2}\partial_{U^n_{\ga_1 \ga_2}}f(v^1,v^2)d\ga_1d\ga_2\Big]\\
&=\sum_{n=1,2}\big<\lambda^n, \partial_{U^n}f(v^1,v^2)\big>_{\mathcal{G}^\omega_0\otimes \mathcal{G}^\omega_0} =Df(v^1,v^2)(\lambda),
\end{split}\]
which is the differential of $f$.  Moreover, $\Gamma^{c}(\omega)$, $\Gamma^{s}(\omega)$, $\Gamma^{1}(\omega)$, $\Lambda^{c}(\omega)$, and $\Lambda^{s}(\omega)$ are defined as follows:
$\forall (j,l)\in\big\{1,\dots,\N{} \big\}^{2}$ and $j\not=l$
\begin{equation*}\begin{split}
\Gamma^{c}_{jl}(\omega)&= \frac{k^4(\omega)}{2\Bh{j}{}\Bh{l}{}}\int_{0}^{+\infty}\mathbb{E}\big[C^\omega_{jl}(0)C^\omega_{jl}(z)\big]\cos\big((\Bh{l}{}-\Bh{j}{})z\big)dz,\\
\Gamma^{c}_{jj}(\omega)&=-\sum_{\substack{l=1\\l\not=j}}^{\N{}}\Gamma^{c}_{jl}(\omega),\\
\Gamma^{s}_{jl}(\omega)&= \frac{k^4(\omega)}{2\Bh{j}{}\Bh{l}{}}\int_{0}^{+\infty}\mathbb{E}\big[C^\omega_{jl}(0)C^\omega_{jl}(z)\big]\sin\big((\Bh{l}{}-\Bh{j}{})z\big)dz, \\
\Gamma^{s}_{jj}(\omega)&=-\sum_{\substack{l=1\\ l\not= j}}^{\N{}}\Gamma^{s}_{jl}(\omega),
\end{split}\end{equation*}
and 
$\forall (j,l)\in\big\{1,\dots,\N{} \big\}^{2}$,
\[\begin{split}
\Gamma^{1}_{jl}(\omega)&= \frac{k^4(\omega)}{2\Bh{j}{}\Bh{l}{}}\int_{0}^{+\infty}\mathbb{E}\big[C^\omega_{jj}(0)C^\omega_{ll}(z)\big]dz,\\
\Lambda^{c}_{j}(\omega)&= \int_{\xi}^{\ko{}}\frac{k^4 (\omega)}{4\sgap \Bh{j}{}}\int_{0}^{+\infty}\mathbb{E}\big[C^\omega_{j\ga'}(0)C^\omega_{j\ga'}(z)\big]\cos\big((\sgap-\Bh{j}{})z\big)dzd\ga', \\
\Lambda^{s}_{j}(\omega)&= \int_{\xi}^{\ko{}} \frac{k^4 (\omega)}{4\sgap \Bh{j}{}}\int_{0}^{+\infty}\mathbb{E}\big[C^\omega_{j\ga'}(0)C^\omega_{j\ga'}(z)\big]\sin\big((\sgap-\Bh{j}{})z\big)dzd\ga',
\end{split}
\]
where the coefficients $C^\omega(z)$ are defined by \eqref{coefcouplVP2}.

From Theorems \ref{thasympP2p1}, we have the following proposition about the autocorrelation function of the transfer operator for the two steps of the time-reversal experiment.
\begin{prop}\label{trpropP2p}
$\forall (y^1,y^2)\in\esp_0 \times \esp_0$ and $\forall \lambda \in \esp_0 \times \esp_0$, the autocorrelation function of the transfer operator for the two steps of the time-reversal experiment as $\e\to 0$ and $\xi\to 0$ is given by
\[\begin{split}
\lim_{\xi\to 0}\lim_{\e\to 0}\mathbb{E}\Big[\big<\emph{\textbf{U}}^{\xi,\e}&(\omega,L)(y^1,y^2),\lambda\big>_{\esp_\xi \otimes \esp_\xi}\Big]=\mathbb{E}\Big[\big<\emph{\textbf{U}}^{0}(\omega,L)(y^1,y^2),\lambda\big>_{\esp_0 \otimes \esp_0}\Big]\\
&=\sum_{j,l=1}^{\N{}} \mathcal{T}^{l}_j(\omega,L)\overline{y^1_l}y^2_l\overline{\lambda_{jj}}+\sum_{\substack{j,m=1\\j\not=m}}^{\N{}}e^{Q_{jm}(\omega)L}\overline{y^1_j}y^2_m\overline{\lambda_{jm}}\\
&\quad+\sum_{j=1}^{\N{}}\int_0^{\ko{}}e^{\frac{1}{2}(\Gamma^c_{jj}(\omega)-\Gamma^1_{jj}(\omega)-\Lambda^{c}_j(\omega))L-\frac{i}{2}(\Gamma^s_{jj}(\omega)-\Lambda^{s}_j(\omega))L}\overline{y^1_j}y^2_{\ga'}\overline{\lambda_{j\ga'}}d\ga'\\
&\quad+\int_0^{\ko{}}\sum_{m=1}^{\N{}}e^{\frac{1}{2}(\Gamma^c_{mm}(\omega)-\Gamma^1_{mm}(\omega)-\Lambda^{c}_m(\omega))L+\frac{i}{2}(\Gamma^s_{mm}(\omega)-\Lambda^{s}_m(\omega))L}\overline{y^1_\ga}y^2_{m}\overline{\lambda_{\ga m}}d\ga\\
&\quad+\int_0^{\ko{}}\int_0^{\ko{}}\overline{y^1_\ga}y^2_{\ga'}\overline{\lambda_{\ga\ga'}}d\ga d\ga'.
\end{split}\]
Here,
\[\begin{split}
Q_{jm}(\omega)=&\frac{1}{2}\big[\Gamma^c_{jj}(\omega)+\Gamma^c_{mm}(\omega)-(\Gamma^1_{jj}(\omega)+\Gamma^1_{mm}(\omega)-2\Gamma^1 _{jl}(\omega))-(\Lambda^{c}_j(\omega)+\Lambda^{c}_m(\omega))\big]\\
&+\frac{i}{2}\big[\Gamma^s_{mm}(\omega)-\Gamma^s_{jj}(\omega)-(\Lambda^{s}_m(\omega)-\Lambda_{l}^{s} (\omega))\big].
\end{split}\]
with $\mathcal{T}^{l}_j(\omega,z)$ is the solution of the coupled power equations
\begin{equation}\label{cpetrceP2p}
\dz\mathcal{T}^{l}_j(\omega,z)=-\Lambda^{c,\xi}_{j}(\omega)\mathcal{T}^{\xi,l}_j(\omega,z)+\sum_{n=1}^{\N{}}\Gamma^c_{nj}(\omega)\big(\mathcal{T}^{l}_n(\omega,z)-\mathcal{T}^{l}_j(\omega,z)\big)
\end{equation}
and $\mathcal{T}^{l}_j(\omega,0)=\delta_{jl}$.
\end{prop}
Let us note that,
\begin{equation}\label{asymptcovtrP2p}\mathcal{T}_j^l(\omega_0,L)=\lim_{\xi\to0}\lim_{\e\to0}\mathbb{E}\Big[\lvert \textbf{T}_j^{\xi,\e}(\omega_0,L)(y^l) \rvert ^2\Big],\end{equation}
with $y^l_j=\delta_{jl}$ and $y^l_\ga=0$ for $\ga \in(0,\ko{})$, is the asymptotic mean mode power of the $j$th propagating mode of the transfer operators at distance $z=L$. The initial condition $y^l$ means that an impulse equal to one charges only the $l$th propagating mode at $z=0$. Equation \eqref{cpetrceP2p} describes the transfer of energy between the propagating and the radiating modes through the  energy transport matrix $\Gamma^{c}(\omega)$ the dissipation coefficients $\Lambda^{c}(\omega)$. These dissipation coefficients resulting from the coupling between the propagating and the radiating modes are responsible to the radiative loss of energy of the propagating modes into the ocean bottom. As we will see in the following section the deterioration of the refocusing is due to these radiative losses which are caused by the random perturbations of the propagation medium.

\subsection{Refocusing in a Random Waveguide}\label{rfcrwintP2p}

Using the asymptotic analysis developed in the previous section, now we are able to describe the asymptotic mean refocused wave. However, let us note that we only need to know the asymptotic mean refocused wave since the refocused wave is self-averaging, which means the refocused wave converges in probability to its asymptotic mean value. We refer to \cite{these} for instance for a complete proof of the self-averaging property. Let us note that the self-averaging property of the time-reversal process has already been observed in many context \cite{papanicolaou3,papanicolaou2,book, papa,gomez}.

Using the change of variable $\omega=\omega_0+\se h$, the refocused wave is given by
\[\begin{split}
p_{TR}\Big(&\frac{t_{obs}}{\e}+\frac{t}{\se},x,L_S\Big) e^{i\omega_0\big(\frac{t_{obs}-t_1}{\e}+\frac{t}{\se}\big)}=\frac{1}{8\pi}\int \overline{\widehat{f}(h)}e^{ih \big(\frac{t_1-t_{obs}}{\se}-t\big)}\\
&\times\big<\textbf{U}^{\xi,\e}(\omega_0+\se h,L)\big(\tilde{a}(\omega_0+\se h),\tilde{b}_x(\omega_0+\se h)\big),\lambda^\e(\omega_0+\se h)\big>_{\mathcal{H}^{\omega_0+\se h}_\xi \otimes \mathcal{H}^{\omega_0+\se h}_\xi}dh.
\end{split}\]
Using Proposition \ref{trpropP2p} we obtain 
\begin{equation}\label{TRfocused}\begin{split}
\mathbb{E}&\Big[p_{TR}\Big(\frac{t_{obs}}{\e}+\frac{t}{\se},x,L_S\Big)\Big] e^{i\omega_0\big(\frac{t_{obs}-t_1}{\e}+\frac{t}{\se}\big)}=\frac{1}{4}\sum_{j,m=1}^{\N{_0}}\sqrt{\frac{\Bh{m}{_0}}{\Bh{j}{_0}}}e^{i(\Bh{m}{_0}-\Bh{j}{_0})\frac{L}{\e}}\\
&\times M_{mj}(\omega_0) K^{\omega_0}_{j,m,L}\ast f \left(\frac{(\beta'_m(\omega_0)-\beta'_j(\omega_0))L+t_1-t_{obs}}{\se}-t\right) \mathbb{E}\Big[\textbf{U}^{\xi}_{jm}(\omega_0,L)\big(\tilde{a}(\omega_0),\tilde{b}_x(\omega_0)\big)\Big]\\
&+\mathcal{O}(\se),
\end{split}\end{equation}
where $K^{\omega_0}_{j,m,L}$ are defined by \eqref{kernelP2p}. Let us note that there is no radiating part in the expression of the refocused wave because all the radiating components of the refocused wave involve a term of the form
\[
\int_\xi^{\ko{}}\phi_{\ga}(\omega,x)\phi_{\ga}(\omega,y)e^{i\sga \frac{L}{\e}}=\mathcal{O}(\e)
\]
uniformly bounded in $x$ and $y$ on bounded subset of $[0,+\infty)^2$. Moreover, we cannot observe the recompression of the radiating components by time reversal because it holds only on a set with null Lebesgue measure.

From \eqref{TRfocused}, we obtain for $m\not=j$
\[\begin{split}
\lim_{\xi \to 0}\lim_{\e\to 0}\mathbb{E}\Big[p_{TR}\Big(\frac{t_{jm}}{\e}&+\frac{t}{\se},x,L_S\Big)\Big]e^{i\omega_0\big(\frac{t_{jm}-t_1}{\e}+\frac{t}{\se}\big)}e^{-i(\Bh{m}{_0}-\Bh{j}{_0})\big(-L_S+\frac{L}{\e}\big)}\\
&=e^{Q_{jm}(\omega_0)L} M_{jm}(\omega_0)K_{j,m,L}^{\omega_0}\ast f(-t) \cdot \frac{1}{4}\phi_{j}(\omega_0,x_0)\phi_m(\omega_0 ,x),  
\end{split}\] 
where the times $t_{jm}$ are defined by \eqref{timetrP2p}. Then, at all these times we can observe the shape of the refocused waves obtained in an homogeneous medium, but with the damping terms $e^{Q_{jm}(\omega_0)L}$. As a result, the amplitude of the coherent refocused waves at times $t_{jm}$ decays exponentially with respect to the propagation distance $L$, and therefore becomes negligible for long propagation distance. More precisely, as we will see in what follows, even for $t_{obs}=t_1$ the amplitude of the refocused wave will decay exponentially fast, but in this case the decay rate is smaller that in the case $t_{obs}=t_{jm}$ ($m\not=j$).

Now, for $t_{obs}=t_1$, we have from Proposition \ref{trpropP2p} a contribution of all the propagating modes
\begin{equation}\label{TRrefocprof}
\lim_{\xi \to 0}\lim_{\e\to 0}\mathbb{E}\Big[p_{TR}\Big(\frac{t_{1}}{\e}+\frac{t}{\se},x,L_S\Big)\Big]e^{i\omega _0\frac{t}{\se}}=f(-t)\cdot \frac{1}{4} \sum_{j,l=1}^{\N{_0}} M_{jj}(\omega_0)\mathcal{T}_j^l(\omega_0,L)\phi_l (\omega_0,x_0)\phi_l(\omega_0,x),\end{equation}
where $\mathcal{T}_j^l(\omega_0,L)$ are the asymptotic mean mode powers \eqref{asymptcovtrP2p} satisfying the coupled power equations \eqref{cpetrceP2p}, and $M_{jj}(\omega_0)$ is defined by \eqref{matrixmirror}.
From this last expression one can see that the  refocusing of time-reversed wave is closely related to the mean mode powers propagation described through the coupled power equations \eqref{cpetrceP2p}. We will see in the following sections that the effects of the radiative losses into the ocean bottom, and caused by the random inhomogeneities, affect the refocused wave in two ways. First, the amplitude of the refocused wave decay exponentially with the propagation distance. This exponential decay of the propagating mode powers is described in \cite{papanicolaou,gomez2}. The second effect is the loss of resolution of the transverse profile. Let us recall that for a waveguide model with a bounded cross section the energy is preserved, so that one can observe the classical time-reversal superresolution effect mainly through a side-lobe suppression produced by the time-reversal mirror \cite{book, papa}.

\subsection{Exponential Decay of the Refocused Wave Amplitude}\label{expodecsection}

The exponential decay rate of the asymptotic propagating mean mode power \eqref{asymptcovtrP2p} is carried out in \cite{gomez2}, but let us recall the result after introducing some notations. 

Let us consider
\[  \mathcal{S}^{\N{_0}}_{+}\!=\!\left\{ X\in\mathbb{R}^{\N{_0}}, \, X_j\geq 0\quad\forall j\in\{1,\dots,\N{_0}\}  \text{ and }\|X\|^2_{2,\mathbb{R}^{\N{_0}}} =\big<X,X\big>_{\mathbb{R}^{\N{_0}}}=1 \right\}\]
with $\big< X,Y\big>_{\mathbb{R}^{\N{_0}}}=\sum_{j=1}^{\N{_0}}X_j Y_j$ for $(X,Y)\in (\mathbb{R}^{\N{_0}})^2$,
and 
\begin{equation}\label{matriceDP2p}
\Lambda^c _d(\omega_0)=diag\big(\Lambda^c (\omega_0),\dots,\Lambda^c _{\N{}}(\omega_0)\big).
\end{equation}
 
\begin{thm} Let us assume that the energy transport matrix $\Gamma^{c}(\omega_0)$ is irreducible. Then, we have
 \[ \lim_{L\to+\infty}\frac{1}{L}\ln\left[ \sum_{j=1}^{\N{}}\mathcal{T}_{j}^l(\omega_0,L) \right]=-\Lambda_{\infty}(\omega_0)\]
with 
\begin{equation*}\Lambda_{\infty}(\omega_0)  =\inf_{X\in \mathcal{S}^{\N{_0}}_{+}} \big< \big(-\Gamma^{c}(\omega_0)+\Lambda^c _d(\omega_0)\big)X , X \big>_{\mathbb{R}^{\N{_0}}}>0.
 \end{equation*}
\end{thm}    
Unfortunately, $\Lambda_{\infty}(\omega_0)$ is not easy to compute, but however we have the following inequalities
\begin{equation}\label{lambdabar} 
0< \Lambda_{min}(\omega_0)=\min_{j\in\{1,\dots,\N{_0}\}}\Lambda^c _j (\omega_0)\quad \leq\quad \Lambda_{\infty}(\omega_0)\quad \leq \quad  \overline{\Lambda}(\omega_0)=\frac{1}{\N{_0}}\sum_{j=1}^{\N{_0}}\Lambda^c_j(\omega_0).
\end{equation}

Let us give two examples for which $\Lambda_{\infty}(\omega_0)$ can be computed explicitly. First, if we assume that the energy transport matrix $\Gamma^{c}(\omega_0)$ can be replaced by $\frac{1}{\tau}\Gamma^{c}(\omega_0)$ with $\tau\ll1$, that is we assume that the mode coupling is strong. In this case one can show \cite{gomez2}, using a probabilistic representation of $\mathcal{T}_{j}^l(\omega_0,L)$ in terms of a jump Markov process, that   
\[\lim_{\tau\to 0}\Lambda^{\tau}_{\infty}(\omega_0)=\overline{\Lambda}(\omega_0),\]
and
\[\lim_{\tau\to 0}\mathcal{T}^{\tau,l}_j (\omega_0,L)=\frac{1}{\N{}}\exp\Big(-\overline{\Lambda}(\omega_0)L\Big).\]
The idea is that we have a strong mixing so that the decays rate averages. On the other hand, if we assume that the energy transport matrix $\Gamma^{c}(\omega_0)$ can be replaced by $\tau\Gamma^{c}(\omega_0)$ with $\tau\ll1$, that is we have a weak mode coupling. In this case, one can show \cite{gomez2} that   
\[\lim_{\tau\to 0}\Lambda^{\tau}_{\infty}(\omega_0)= \Lambda_{min} (\omega_0),\]
The idea is that the modes coupling is too weak to provide any mixing so that the decay rate for the $j$th mode is given by $\Lambda^c_j(\omega_0)$.

As we will see in what follows, the loss of resolution in the time-reversal experiment does not depend only on the mode coupling between the radiative and the propagating modes, but also on the energy transfer between the propagating modes, and described by the transfer matrix $\Gamma^{c}(\omega_0)$.

\subsection{Refocused Transverse Profile without Radiative Losses}\label{sectnorad}

Before, studying the effect of the radiative losses on the focusing quality, let us recall the basic result in the case of negligible coupling between the propagating and the radiating modes \cite{book,papa}. 

In this section and the following one let us note
\begin{equation}\label{transprof}H^{\alpha_M}_{x_0}(\omega_0,x,L)=\frac{1}{4}\sum_{j,l=1}^{\N{_0}}M_{jj}(\omega_0)\mathcal{T}^{l}_j(\omega_0,L)\phi_l(\omega_0,x_0)\phi_l(\omega_0,x),\end{equation}
the transverse profile of \eqref{TRrefocprof}, where $M_{jj}(\omega_0)$ is defined by \eqref{matrixmirror}, and $\alpha_M$ represents the order of magnitude of the time reversal mirror. 

In this section, we assume that the radiative losses in the ocean bottom caused by the random perturbations of the propagation medium are negligible, that is $\Lambda^c(\omega_0)=0$. Consequently, the asymptotic mean mode powers satisfy \eqref{cpetrceP2p} without the radiating coefficient  $\Lambda^c(\omega_0)$. Moreover, according to \cite{book,papa}, using a probabilistic representation of the propagating mean mode powers in terms of a jump Markov, we have
\[\underset{L\to +\infty}{\lim}\mathcal{T}^{l}_j(\omega_0,L)=\frac{1}{N(\omega_0)}.\]
This result describes the asymptotic behavior of the propagating mode powers for long propagation distances. In this asymptotic the energy is uniformly distributed over all the propagating modes and there is no loss of energy, in other words there is no loss of information. As a result, the refocused transverse profile \eqref{transprof} is given by
\begin{equation}\label{profnorad}\underset{L\to +\infty}{\lim} H^{\alpha_M}_{x_0}(\omega_0,x,L)=\frac{1}{4}\sum_{j=1}^{\N{_0}}M_{jj}(\omega_0)\frac{1}{N}\sum_{l=1}^{\N{_0}}\phi_l(\omega_0,x_0)\phi_l(\omega_0,x),\end{equation}
which is closed, up to a multiplicative constant, to the sinc profile according to Proposition \ref{prop17P2p}. However, let us note that in contrast to Section \ref{rfhwintP2p} for a homogeneous waveguide the resolution does not depend on the time-reversal mirror. Let us also note that this transverse profile corresponds to the projection of $\delta(x-x_0)$ over the propagating modes, which corresponds to a point source localized at $x=x_0$. In absence of energy loss, or information loss, the refocused transverse profile, which is closed to the sinc profile, is then the best profile which can be expected for the time-reversal experiment in a random waveguide. However, to observe better refocusing properties, the waveguide setup has to be modified as described in \cite{gomez}.

\subsection{Example of a Refocused Transverse Profile without Loss of resolution}

In this section, we show that the loss of resolution is not only due to the mode coupling between the radiative and the propagating modes. In fact, in the first example introduced in Section \ref{expodecsection} the resolution quality is not affected, we still obtain the profile \eqref{profnorad} but with a damping term decaying exponentially fast, and with a decay rate depending on the mode  coupling between the radiative and the propagating modes.

If the propagating mode coupling process is stronger than the radiative losses, that is we consider $\frac{1}{\tau}\Gamma^{c}(\omega_0)$ with $\tau\ll1$,  we have  
\[\lim_{\tau\to 0}\mathcal{T}^{\tau,l}_j (\omega_0,L)=\frac{1}{\N{}}\exp\Big(-\overline{\Lambda}(\omega_0)L\Big),\]
where $\overline{\Lambda}(\omega_0)$ is defined in \eqref{lambdabar}. Consequently, the refocused transverse profile \eqref{transprof} is given by
\[\underset{L\to +\infty}{\lim} H^{\alpha_M}_{x_0}(\omega_0,x,L)=\frac{1}{4}\exp\Big(-\overline{\Lambda}(\omega_0)L\Big)\sum_{j=1}^{\N{_0}}M_{jj}(\omega_0)\frac{1}{\N{_0}}\sum_{l=1}^{\N{_0}}\phi_l(\omega_0,x_0)\phi_l(\omega_0,x),\]
which is the refocused transverse profile obtained without any radiative losses \eqref{profnorad} but with a damping term given by the average radiative losses. Because of the strong mixing property the energy carried by the propagating mode is uniformly distributed as  discussed in Section \eqref{expodecsection}, but decaying exponentially fast, so that the transverse profile is the projection of $\delta(x-x_0)$ over the propagating modes.

\subsection{ Loss of Resolution for the Refocused Transverse Profile in the Continuum Limit}

This section is devoted to the study of the loss of resolution of the refocused transverse profile. Considering, a nearest neighbor coupling mechanism in \eqref{cpetrceP2p}, we use the continuum approximation ($\N{_0}\gg1$) of this equation to describe the refocused transverse profile.

\subsubsection{Nearest Neighbor Mode Coupling}

In the following sections we investigate the loss of resolution of the refocused wave. To lighten the loss of resolution and in order to give simple representation of refocused the transverse profile \eqref{transprof}, we consider the simplified coupled power equation
\begin{equation}\label{simpcpe}\begin{split}
\dz \mathcal{T}_{N} ^l (\omega_0,z)&=-\Lambda^c _{N}(\omega_0)\mathcal{T}_{N} ^l (\omega_0,z)+\Gamma^c _{N-1\,N}(\omega_0)\left(\mathcal{T}_{N-1} ^l (\omega_0,z)-\mathcal{T}_{N} ^l (\omega_0,z)\right),\\
\dz \mathcal{T}_{j} ^l (\omega_0,z)&=\Gamma^c _{j-1\,j}(\omega_0)\left(\mathcal{T}_{j-1} ^l (\omega_0,z)-\mathcal{T}_{j} ^l (\omega_0,z)\right)\\
&+\Gamma^c _{j+1\,j}(\omega_0)\left(\mathcal{T}_{j+1} ^l (\omega_0,z)-\mathcal{T}_{j} ^l (\omega_0,z)\right) \text{ for }j\in\{2,\dots,N-1\},\\
\dz \mathcal{T}_{1} ^l (\omega_0,z)&=\Gamma^c _{2\,1}(\omega_0)\left(\mathcal{T}_{2} ^l (\omega_0,z)-\mathcal{T}_{1} ^l (\omega_0,z)\right),
\end{split}\end{equation} 
with $\mathcal{T}_{j} ^l (\omega_0,0)=\delta_{jl}$,
that is we only consider  a nearest neighbor mode coupling. As a result, only the highest propagating mode can  still be coupled with the radiating modes since it is the closest mode to the continuous spectrum of the Pekeris operator. This simplified version of the coupled power equation can be rigorously derived using a band-limiting idealization \cite{gomez2,papanicolaou}, that is the power spectral density of the random perturbations in the transverse direction has a compact support. This simplification will be the starting point to study the refocusing properties using the continuous diffusion approximations of \eqref{simpcpe} introduced in \cite{gomez2}. These continuum approximations allow to exhibit in a very simple way the effects of the radiative losses on the transverse wave refocusing. First, let us recall what happened if the radiation losses are negligible, that is we neglect the coupling between the propagating and the radiating modes. Afterward, we describe the influence of the radiative losses in the ocean bottom on the time-reversal refocusing.

\subsubsection{Refocused Wave in the Continuum Limit with Negligible Radiation Losses}\label{mrfnrlintP2p}

In this section, we study the transverse profile of the refocused wave in the case where the radiation losses are negligible, that is $\Lambda^{c}_N(\omega)=0$. In this case, the two following propositions describe the refocusing properties of \eqref{transprof}.

According to Proposition \ref{prop17P2p}, let us recall that the size of the focal spot is of order the carrier wavelength of the ocean $\lambda_{oc}=2\pi c/(n_1 \omega_0)$, which tends to $0$ in this continuum limit $\N{_0}\gg1$. Let us also recall that the continuum limit $\N{_0}\gg1$ corresponds to the case $\omega_0\nearrow +\infty$. Consequently, we study the transverse profile \eqref{transprof} of the refocused wave in a window of size $\lambda_{oc}$ centered around $x_0$.

\begin{prop}\label{propsinctrP2p}
For $\alpha_M\in[0,1)$, with negligible radiation losses, the transverse profile of the refocused wave in the continuum limit $\N{_0}\gg 1$ is given by
\[
\lim_{\omega_0 \to +\infty}\frac{\lambda_{oc}^{1-\alpha_M}}{\theta}H^{\alpha_M}_{x_0}\Big(\omega_0,x_0+\frac{\lambda_{oc}}{\theta}\tilde{x},L\Big)=\frac{\tilde{d}_2+\tilde{d}_1}{d}\emph{sinc}(2\pi\tilde{x}).
\]
\end{prop}
\begin{figure}\begin{center}
\includegraphics*[scale=0.6]{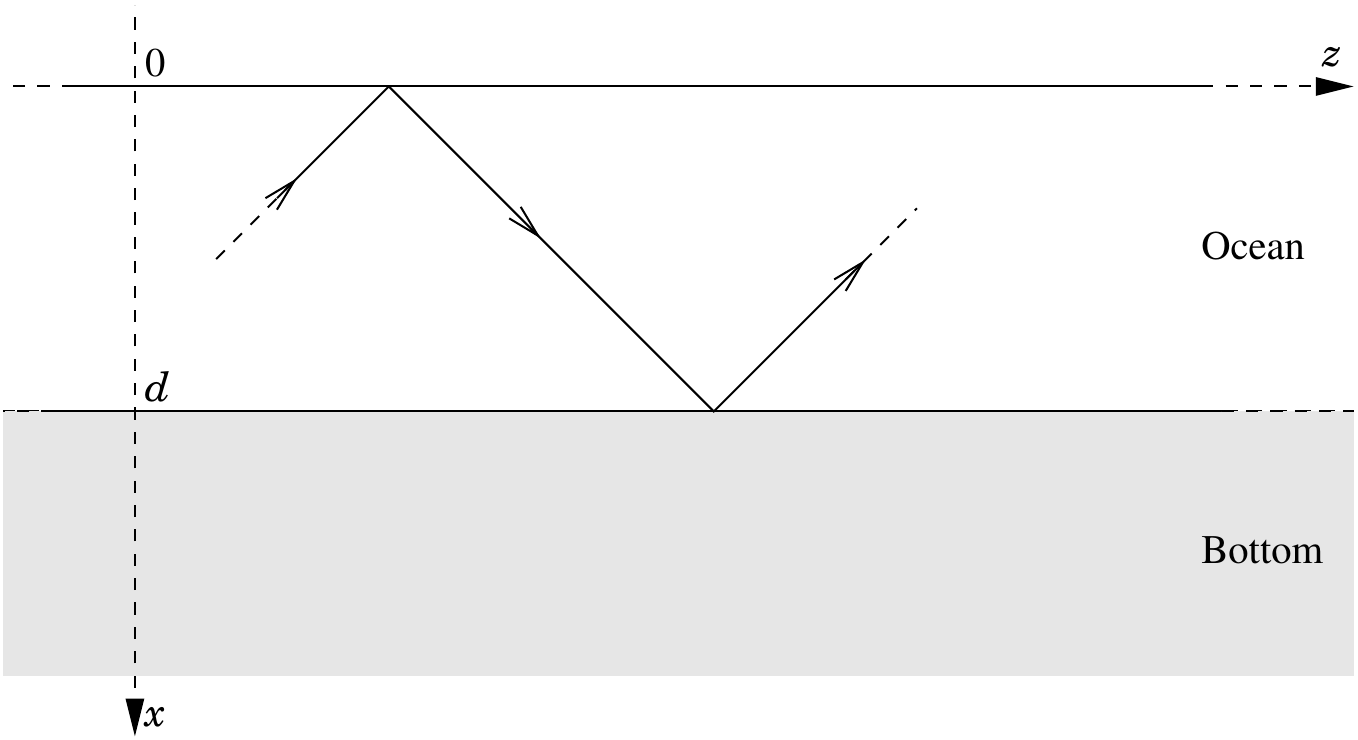}
\end{center}
\caption{\label{perterad2P2p} Illustration of negligible radiation losses in the shallow-water random waveguide model.}
\end{figure}
The transverse profile of the refocused wave is studied using the continuous diffusion approximation model of the coupled power equation \eqref{simpcpe} introduced in \cite[Theorem 6.5]{gomez2}. Without radiative losses the energy is conserved (see Figure \ref{perterad2P2p}) so that we obtain the best transverse profile, which is the sinc profile. In the same way, we have the following result for $\alpha_M=1$, that is the carrier wavelength and the magnitude of the time-reversal mirror are of same order. In order to study this case, let us introduce some notations. Let $\mathcal{E}=\bigcup_{M\geq1}\mathcal{E}_M$, where
\[\mathcal{E}_{M}=\left\{ \sum_{j=1}^M \alpha_j \phi_j,\quad (\alpha_j)_j \in\mathbb{R}^M  \right\},\quad
\text{and}\quad\phi_j(x)=\sqrt{\frac{2}{d}}\sin\Big(j\frac{\pi}{d}x\Big)\quad\forall x\in[0,d],\forall j\geq 1. \]
Let us remark that $(\phi_j)_j$ is a basis of $L^2(0,d)$. 
\begin{prop}\label{transprof20P2p}
For $\alpha_M=1$, in the continuum limit $\N{_0}\gg1$, we have 
\[\lim_{\omega_0\to +\infty}H^{1}_{x_0}(\omega_0,.,L)-\tilde{H}^{1}_{x_0}(\omega_0,.,L)=0\]
in $\mathcal{E}'$, and where
 \[\lim_{\omega_0 \to +\infty}\tilde{H}^{1}_{x_0}\Big(\omega_0,x_0+\frac{\lambda_{oc}}{\theta}\tilde{x},L\Big)=\theta\frac{\tilde{d}_2+\tilde{d}_1}{d}\emph{sinc}(2\pi\tilde{x}).\]
\end{prop}
These results are consistent with the ones obtained Section \ref{sectnorad} and in \cite[Chapter 20]{book} and \cite{papa}, where the authors have obtained the $\textrm{sinc}$ function for transverse profile. The most important fact is that, according to Proposition \ref{prop17P2p} the transverse profile of the refocused wave does not depend on the magnitude of the time-reversal mirror, we obtain the sinc profile $\forall \alpha_M\in[0,1]$.

However, when the radiation losses are not negligible anymore the refocusing properties will be affected in two ways. First, the amplitude of the wave decrease exponentially fast with the size of the random section, and second, the refocusing quality itself is deteriorated by its losses.

\subsubsection{Refocused Wave in the Continuum Limit  with Radiation Losses}\label{mrfrsintP2p}

In addition to the exponential decay of the refocused wave amplitude, the second main effect of the radiative losses is a deterioration of the focusing quality described in the two following theorems. To study the refocused transversed profile \eqref{transprof} we use the continuous diffusion approximation model introduced in \cite[Theorem 6.3]{gomez2} of the coupled power equation \eqref{simpcpe}. According to Proposition \ref{prop17P2p}, let us recall one more time that the size of the focal spot is of order the carrier wavelength of the ocean $\lambda_{oc}$, which tends to $0$ in this continuum limit $\N{_0}\gg1$  ($\omega_0\nearrow +\infty$). Consequently, we study the transverse profile \eqref{transprof} of the refocused wave in a window of size $\lambda_{oc}$ centered around $x_0$.

\begin{thm}\label{transprofP2p}
For $\alpha_M\in[0,1)$, the transverse profile of the refocused wave in the continuum limit $\N{_0}\gg 1$ is given by
\[
\lim_{\omega_0 \to +\infty}\frac{\lambda_{oc}^{1-\alpha_M}}{\theta}H^{\alpha_M}_{x_0}\Big(\omega_0,x_0+\frac{\lambda_{oc}}{\theta}\tilde{x},L\Big)=\frac{\tilde{d}_2+\tilde{d}_1}{d}H(\tilde{x},L),
\]
where 
\begin{equation}\label{refocprofH}H(\tilde{x},L)=\int_0^1 \mathcal{T}_1 (L,u)\cos(2\pi u\tilde{x})du,\end{equation}
and
$\mathcal{T}_1 (L,u)$ is the solution of
\[\frac{\partial}{\partial z} \mathcal{T}_1 (z,u)= \frac{\partial}{\partial u}\left(a_{\infty}(\cdot)\frac{\partial}{\partial u}  \mathcal{T}_1\right)(z,u), \]
with the boundary conditions:
\[ \frac{\partial}{\partial u} \mathcal{T}_1 (z,0)=0,\quad  \mathcal{T}_1 (z,1)=0 \quad \text{and}\quad \mathcal{T}_1 (0,u)=1,\]
 $\forall z>0$. Here,
\[a_{\infty}(u)=\frac{a_0 }{1-\left(1-\frac{\pi^2 }{a^2 d^2} \right)(\theta u)^2},\]
with $a_0  = \frac{ \pi^2  S_{0} }{ 2 a n_1 ^4 d^4 \theta ^2}$, $\theta=\sqrt{1-1/n^2_1}$, $S_0 =\int_0^d\int_0^d \ga_0(x_1,x_2)\cos\big(\frac{\pi}{d}x_1\big)\cos\big(\frac{\pi}{d}x_2\big)dx_1 dx_2$. $n_1$ is the index of refraction in the ocean section $[0,d]$, $1/a=l_{z,x}$ is the correlation length of the random inhomogeneities in the longitudinal direction, and $\ga_0$ is the covariance function of the random inhomogeneities in the transverse direction.
\end{thm}
\begin{figure}\begin{center}
\includegraphics*[scale=0.6]{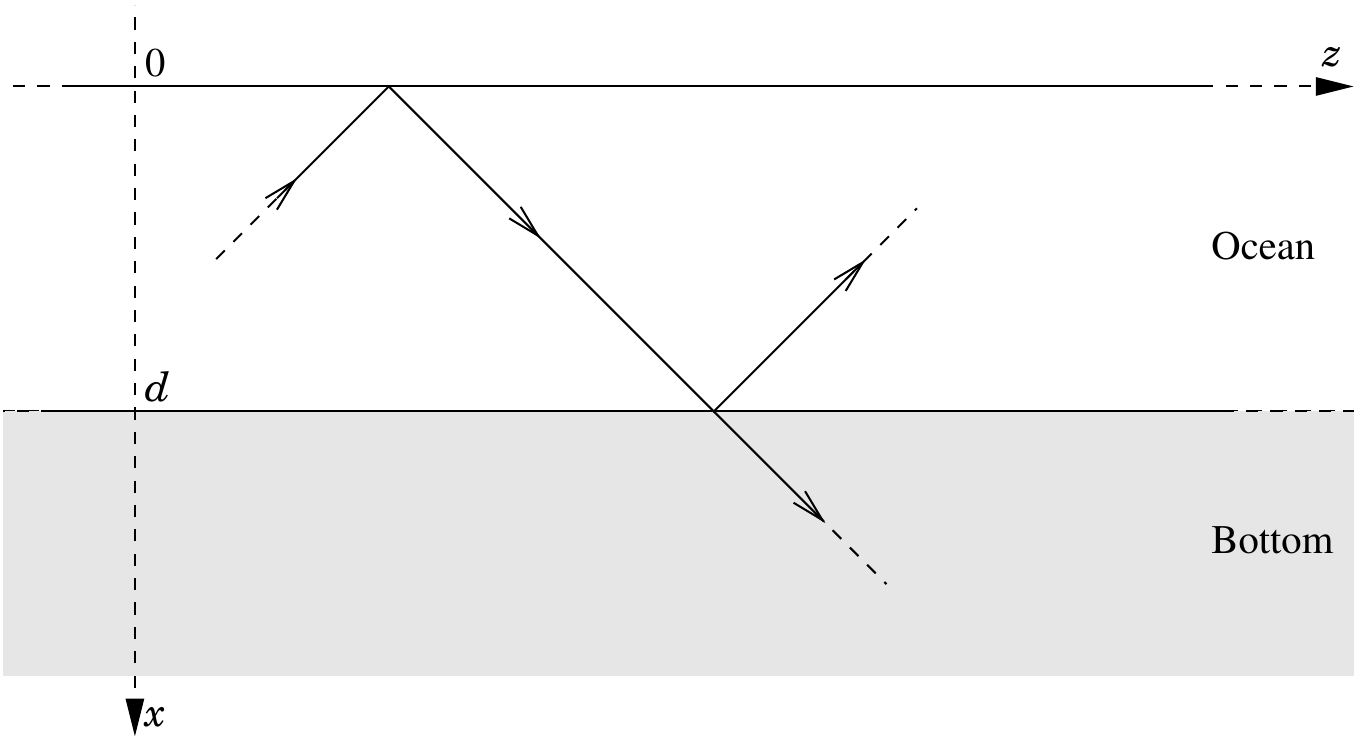}
\end{center}
\caption{\label{perterad1P2p} Illustration of the radiative loss in the shallow-water random waveguide model.}
\end{figure}
The proof of Theorem \ref{transprofP2p} is given in Section \ref{proof2}. Consequently, the transverse profile of the refocused wave can be expressed in terms of the diffusive continuous model, with a reflecting boundary condition at $u=0$ (the top of the waveguide) and an absorbing boundary condition at $u=1$ (the bottom of the waveguide) which represents the radiative loss (see Figure \ref{perterad1P2p}).  

Thanks to Theorem \ref{transprofP2p}, we can now give a simple explanation of this refocusing degradation. First, let us recall that the quality of the refocusing depends much more on the high propagating modes that the lower ones.  As illustrated in Figure \ref{figproftr12P2p}, the absorbing boundary condition at  $u=1$, describing the radiation losses, degrades the information carried by the high propagating modes and it is getting worst as the size of the random section becomes large. As a result, the function $\mathcal{T}_1(L,u)$ in \eqref{refocprofH} plays the role of a low pass filter, which therefore degrades the refocusing quality as illustrated in Section \ref{sectillnumP2p} (Figure \ref{figproftr12P2p} and Figure \ref{resolLgrdP2p}). Consequently, the radiation losses degrade the quality of the refocusing in two ways: the amplitude of the refocused wave decays exponentially with the propagation distance, and the width of the focal spot increases and converges to an asymptotic value as $L\to +\infty$ that is significantly larger than the diffraction limit $\lambda_{oc}/(2\theta)$, where $\lambda_{oc}$ is the carrier wavelength in the ocean section $[0,d]$ (see Section \ref{sectillnumP2p}). 

Let us note that in the case of negligible radiation losses we have the same diffusive continuous model as in Theorem \ref{transprofP2p} but with two reflecting boundary conditions at $u=0$ and $u=1$, (see Section \ref{proof3} in Appendix to see the full details). However, in this case the solution of the diffusion equation is trivial because of the energy conservation. As a result the main difference between the cases with or without radiation losses is the low pass filter effect produced by the  absorbing boundary condition at $u=1$ (the bottom of the waveguide), which breaks the energy conservation.
  
Let us also remark that the exponential decay of the refocused wave amplitude can also be obtained directly from the diffusive continuous model \cite[Theorem 6.3]{gomez2} described in Theorem \ref{transprofP2p}. Finally, we also have the same result for $\alpha_M=1$, that is the carrier wavelength and the magnitude of the time-reversal mirror are of same order.

\begin{thm}\label{transprof2P2p}
For $\alpha_M=1$, in the continuum limit $\N{_0}\gg1$, we have 
\[\lim_{\omega_0\to +\infty}H^1_{x_0}(\omega_0,.,L)-\tilde{H}^1_{x_0}(\omega_0,.,L)=0\]
in $\mathcal{E}'$, which is the topological dual of $\mathcal{E}$ equipped with the weak topology, and where
 \[
\lim_{\omega_0 \to +\infty}\tilde{H}^{1}_{x_0}\Big(\omega_0,x_0+\frac{\lambda_{oc}}{\theta}\tilde{x},L\Big)=\theta\frac{\tilde{d}_2+\tilde{d}_1}{d}H(\tilde{x},L).
\]
Here, $H(\tilde{x},L)$ is defined in Proposition \ref{transprofP2p}.
\end{thm}

Let us remark that in the case of a random waveguide, the order of magnitude $\alpha_M$ of the time-reversal mirror plays no role in the transverse profile compared to the homogeneous case.

\subsubsection{Numerical Illustrations}\label{sectillnumP2p}

In this section we illustrate the spatial focusing of the refocused wave around the source location. First, we represent the evolution of $\mathcal{T}_1(L,u)$, in presence of radiation losses, with respect to $L$. Here, $\mathcal{T}_1(L,u)$ is the mean mode power for the $[\N{_0}u]$th propagating mode in the continuum limit $\N{_0}\gg1$, which is the solution of the partial differential equation in Theorem \ref{transprofP2p}.

 Second, we represent the transverse profile $H(\tilde{x},L)$ defined by \eqref{refocprofH} of the refocused wave, and finally we illustrate the resolution of the refocused wave as the propagation distance $L$ becomes large. 

In this section, we consider the following values of the parameters. For the sake of simplicity, we take $a_0=1$, and the correlation length of the random inhomogeneities in the longitudinal direction is $1$ ($a=1$). Moreover, we take $n_1=2$ for index of refraction in the ocean section $[0,d]$, and depth $d=20$.     

\begin{figure}
\begin{center}\includegraphics*[scale=0.4]{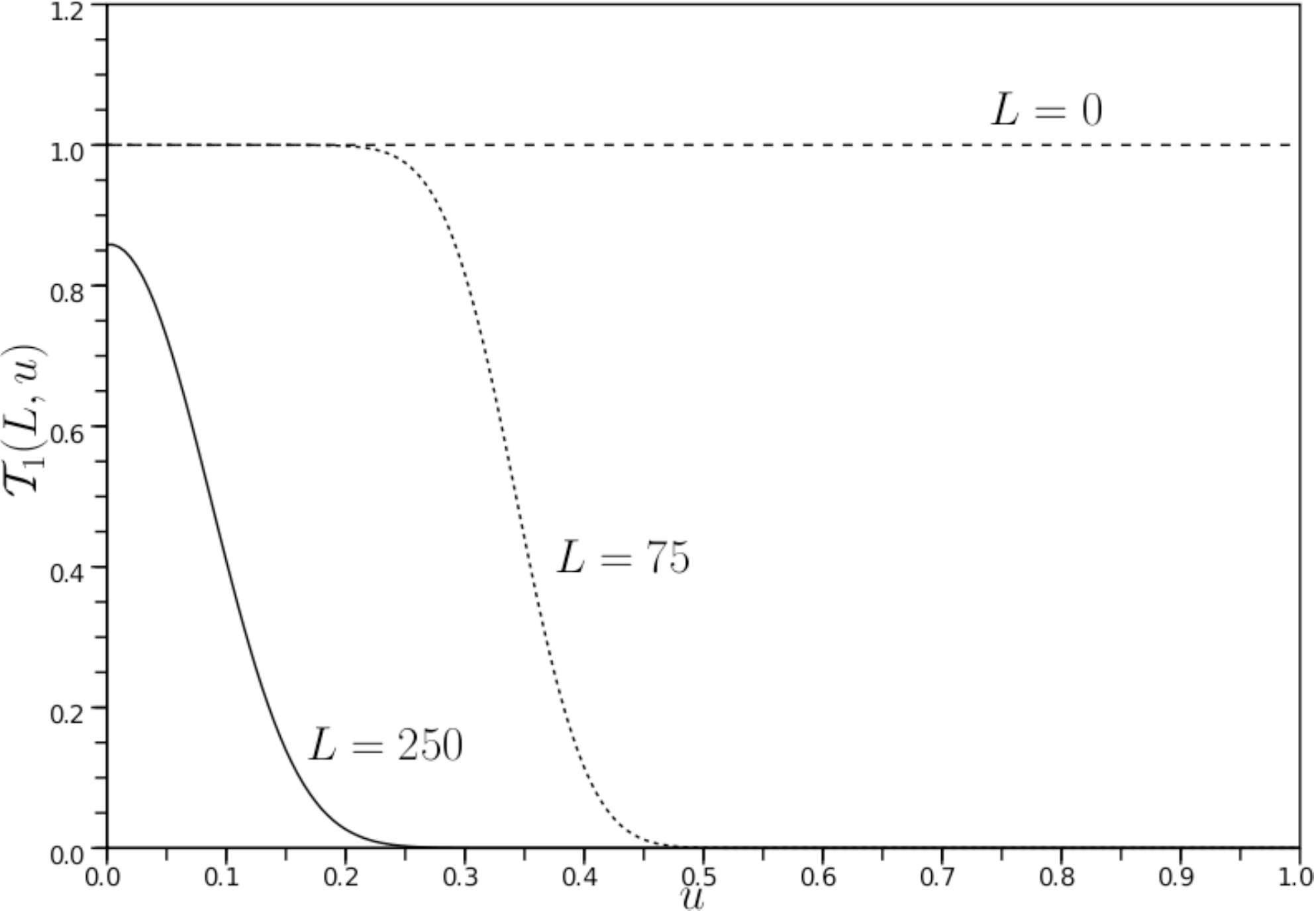}\end{center}
\caption{\label{tau1P2p} Representation of $\mathcal{T}_1(L,u)$, in the presence of radiation losses, with respect to the propagation distance $L$.}
\end{figure}

\begin{figure} \begin{center}
\begin{tabular}{cc}
\includegraphics*[scale=0.37]{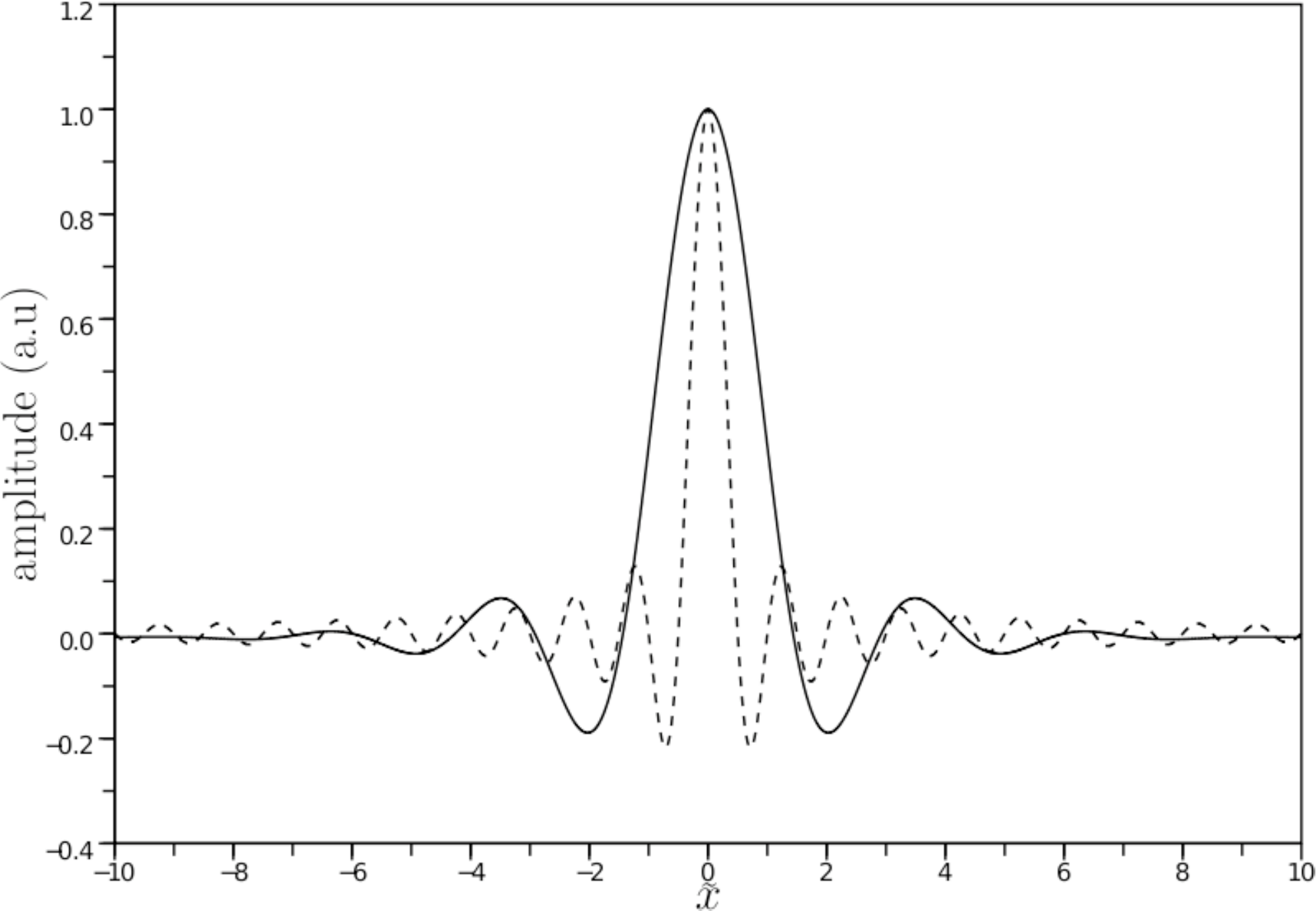}&\includegraphics*[scale=0.4]{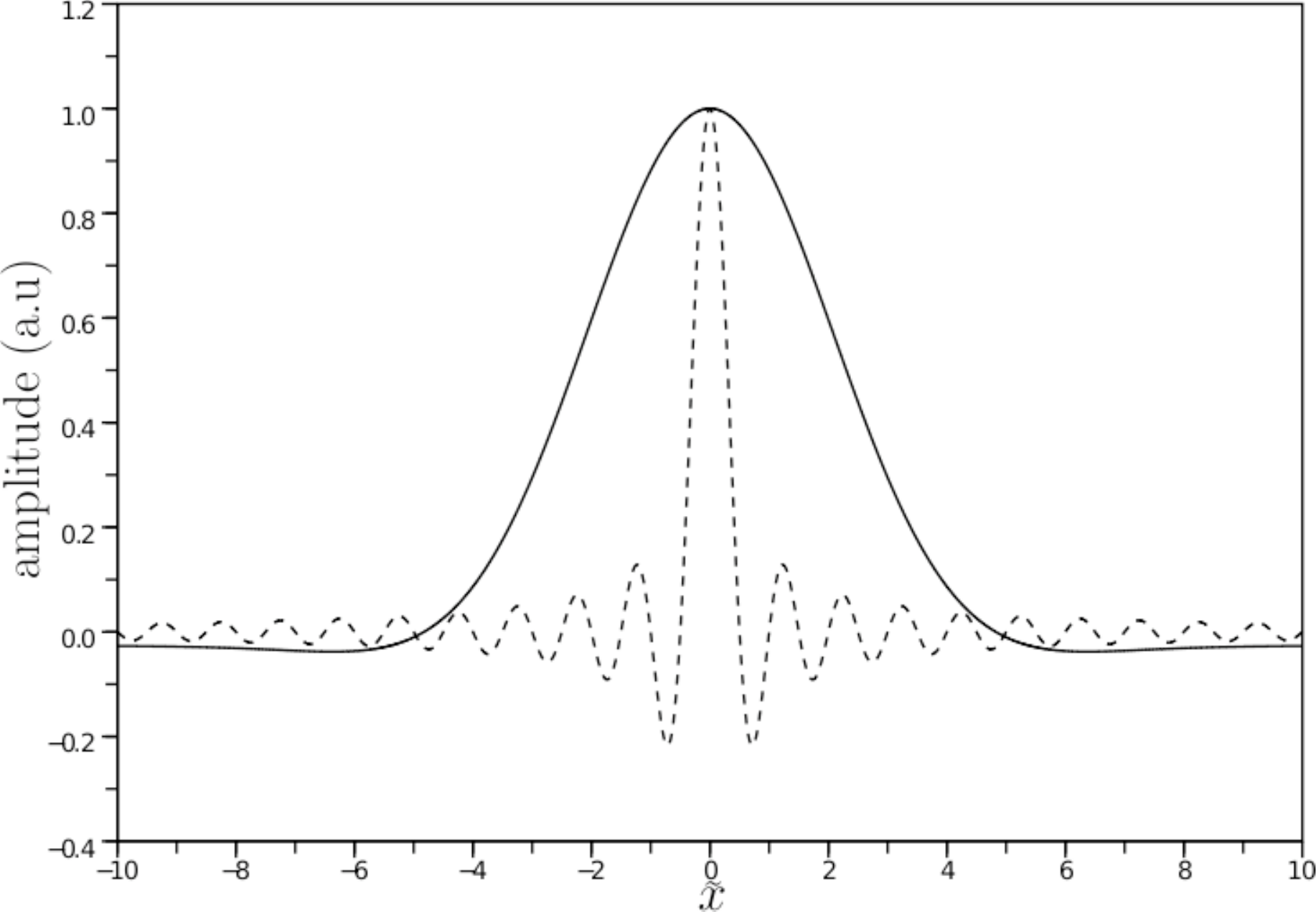}\\
$(a)$&$(b)$
\end{tabular}
\end{center}
\caption{\label{figproftr12P2p}
Normalized transverse profile. In $(a)$ and $(b)$ the dashed curves are the transverse profiles in the case where the radiation losses are negligible, and the solid curves represent the transverse profile $H(\tilde{x},L)$. In $(a)$ we represent $H(\tilde{x},L)$ with $L=75$, and in $(b)$ we represent $H(\tilde{x},L)$ with $L=250$. 
}\end{figure}

\begin{figure}
\begin{center}\includegraphics*[scale=0.4]{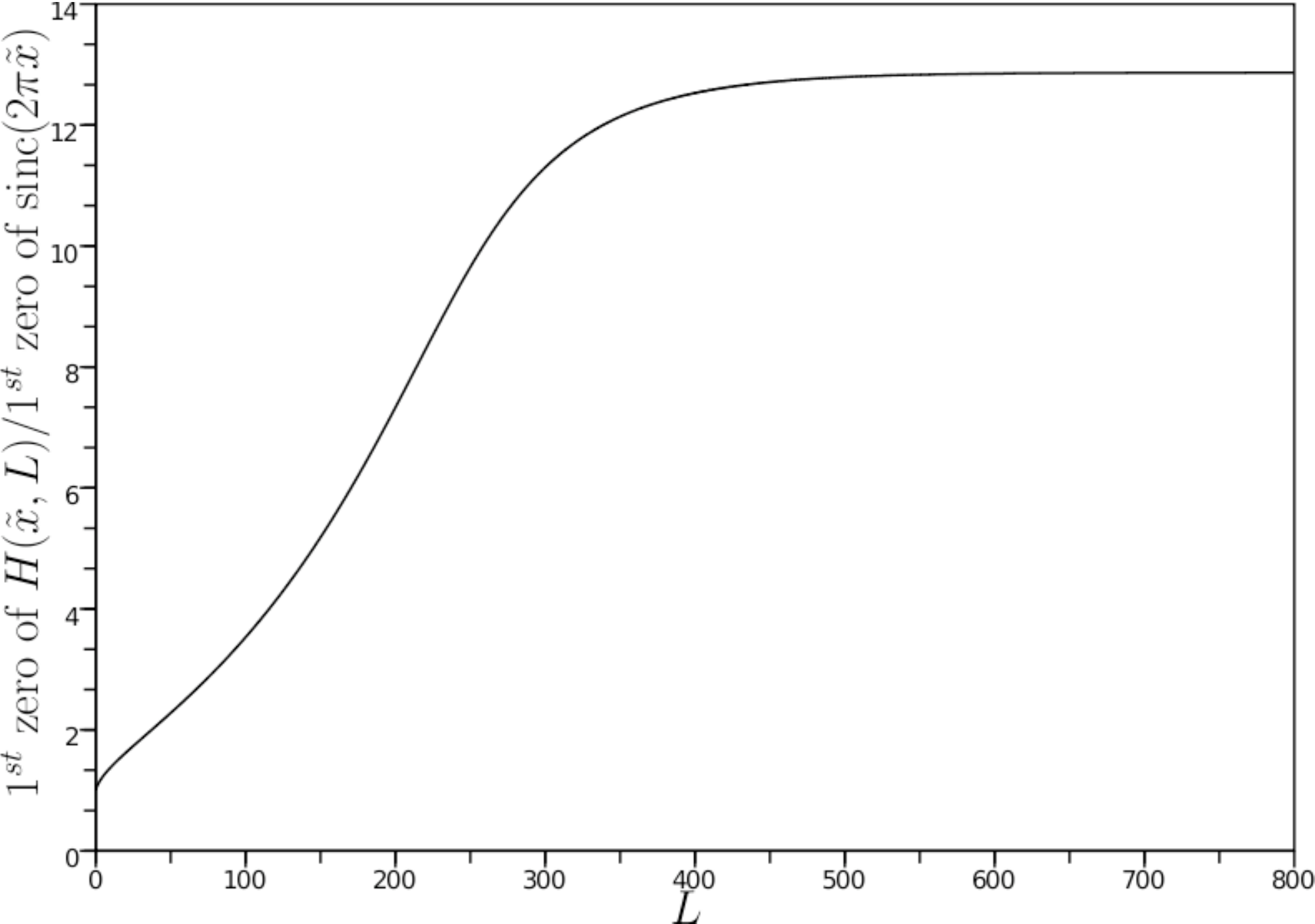}\end{center}
\caption{\label{resolLgrdP2p} Representation of the evolution of the resolution with respect to the propagation distance $L$.}
\end{figure}

We saw in Theorem \ref{transprofP2p} and Theorem \ref{transprof2P2p} that $\mathcal{T}_1(L,u)$, in the presence of radiation losses, plays an important role in the transverse profile of the refocused wave. In Figure \ref{tau1P2p},
 we illustrate the influence of the radiation losses on $\mathcal{T}_1(L,u)$ as the propagation distance $L$ increase. As we can see in Figure \ref{figproftr12P2p} 
and Figure \ref{resolLgrdP2p},
 the radiation losses degrade the quality of the refocusing. Moreover, for $L\gg 1$, one can see a threshold of the quality of the resolution since
\[H_{x_0}(\tilde{x},L)\underset{L\gg1}{\simeq}e^{\lambda_1 L}\int_0^1 \phi_{\infty,1}(v)dv \int_0^1 \phi_{\infty,1}(u)\cos(2\pi\tilde{x}u)du,\]
where $\lambda_1<0$ is a simple eigenvalue, and the largest one, with corresponding eigenvector $\phi_{\infty,1}$ of the diffusion operator introduced in Theorem  \ref{transprofP2p} (see \cite[Lemma 2.2]{these}).

\section*{Conclusion}

In this paper we have analyzed the time reversal of waves of a broadband pulse in an underwater acoustic channel with random perturbations. In this context, using the continuous diffusive models developed in \cite{gomez2}, describing the mode-power coupling between the propagating and radiating modes,  we describe in a simple way the refocused transverse profile (Theorem \ref{transprofP2p} and Theorem \ref{transprof2P2p}) to lighten the negative effects of the radiative losses on the time-reversal refocusing property. We have seen that radiation losses degrade the quality of the refocused transverse profile as the propagation distance increases. First, the amplitude of the refocused wave decays exponentially with the propagation distance (Section \ref{expodecsection}). Second, using a low pass filter representation of the refocused transverse profile (Theorem \ref{transprofP2p}), we have shown that the width of the main focal spot increases and converges to an asymptotic value as the size of the random section increase, which is significantly larger than the diffraction limit $\lambda_{oc}/(2\theta)$ obtained in Proposition \ref{prop17P2p} (where $\lambda_{oc}$ is the carrier wavelength in the ocean section with index of refraction $n_1$, and $\theta=\sqrt{1-1/n_1^2}$).

\section{Appendix}

\subsection{Proof of Proposition \ref{prop17P2p}}\label{proof0}

Let us first note that according to Lemma \ref{coefgP2} and \eqref{approxvp}, we have
\begin{equation}\label{preveq}\sup_{j\in\{1,\dots,N-N^\alpha\}}\lvert A_j^2-\frac{2}{d}\rvert=\mathcal{O}(N^{\alpha-1}) \end{equation}
where $A_j$ is defined by \eqref{coefajP2}. Then, we split the transverse profile in two part in order to use Lemma \ref{coefgP2} and \eqref{preveq},
\[\begin{split}
\frac{d}{\tilde{d}_2+\tilde{d}_1}&\frac{\lambda_{oc}^{1-\alpha_M}}{\theta}H^{\alpha_M}_{x_0}\Big(\omega_0,x_0+\frac{\lambda_{oc}}{\theta}\tilde{x}\Big)=\frac{\lambda_{oc}}{2\theta}\left[\sum_{j=1}^{N-[N^\alpha]}+\sum_{j=N-[N^\alpha]+1}^{N}\right]\phi_j(\omega_0,x_0)\\
&\phi_j(\omega_0,x_0+\lambda_{oc}\tilde{x}/\theta)\frac{A_j^2 d}{2}\Big[1-\cos\Big(\sigma_j\frac{2d_M+\lambda_{oc}^{\alpha_M}(\tilde{d}_2-\tilde{d}_1)}{d}\Big)\textrm{sinc}\Big(\sigma_j\frac{\lambda_{oc}^{\alpha_M}(\tilde{d}_2+\tilde{d}_1)}{d}\Big)\Big].
\end{split}\]
so that the second sum on the right of the previous equality is of order $\mathcal{O}(N^{\alpha-1})$. Moreover,
\[\begin{split}
\Big\vert  \frac{\lambda_{oc}}{2\theta}&\sum_{j=1}^{N-[N^\alpha]}\phi_j(\omega_0,x_0)\phi_j(\omega_0,x_0+\lambda_{oc}\tilde{x}/\theta)\\
&\times \frac{A_j^2 d}{2}\cos\Big(\sigma_j\frac{2d_M+\lambda_{oc}^{\alpha_M}(\tilde{d}_2-\tilde{d}_1)}{d}\Big)\textrm{sinc}\Big(\sigma_j\frac{\lambda_{oc}^{\alpha_M}(\tilde{d}_2+\tilde{d}_1)}{d}\Big) \Big\vert \leq K \lambda_{oc}^{1-\alpha_M}\ln(N).
\end{split}\]
Now, for the first sum of the previous equality we have 
\[\phi_j (\omega_0,x_0)\phi_j\big(\omega_0,x_0+\frac{\lambda_{oc}}{\theta}\tilde{x}\big)=\frac{A^2_j}{2}\Big[\cos\Big(\sigma_j\frac{\lambda_{oc}}{\theta d}\tilde{x}\Big)-\cos\Big(\sigma_j \frac{2x_0 +\lambda_{oc}\tilde{x}/\theta}{d}\Big)\Big],\]
and 
\[\begin{split}
\cos\Big(\sigma_j \frac{2x_0 +\lambda_{oc}\tilde{x}/\theta}{d}\Big)=&\cos\Big((\sigma_j-j\pi) \frac{2x_0 +\lambda_{oc}\tilde{x}/\theta}{d}\Big)\cos\Big(j\pi \frac{2x_0 +\lambda_{oc}\tilde{x}/\theta}{d}\Big)\\
&-\sin\Big((\sigma_j-j\pi) \frac{2x_0 +\lambda_{oc}\tilde{x}/\theta}{d}\Big)\sin\Big(j\pi \frac{2x_0 +\lambda_{oc}\tilde{x}/\theta}{d}\Big).
\end{split}\]
Then, using the Abel transform and Lemma \ref{coefgP2}, we get
\[\lambda_{oc}\Big\lvert \sum_{j=1}^{N-[N^\alpha]}\cos\Big(\sigma_j \frac{2x_0 +\lambda_{oc}\tilde{x}/\theta}{d}\Big) \Big\rvert \leq K N^{\frac{1}{2}-\frac{3}{2}\alpha}.\]
Moreover, using \eqref{preveq}, we also have 
\[
\frac{\lambda_{oc} d}{8\theta}\sum_{j=1}^{N-[N^\alpha]} A^4_j \cos\Big(\sigma_j\frac{\lambda_{oc}}{\theta d}\tilde{x}\Big) =\frac{\lambda_{oc} }{2\theta d}\sum_{j=1}^{N-[N^\alpha]}\cos\Big(2\frac{j}{N}\pi\tilde{x}\Big) + \mathcal{O}(N^{\alpha-1}),
\]
with
\[\lim_{\omega_0\to +\infty}\frac{\lambda_{oc} }{2\theta d}\sum_{j=1}^{N-[N^\alpha]}\cos\Big(2\frac{j}{N}\pi\tilde{x}\Big) =\int_0^1 \cos(2u\pi\tilde{x})du=\textrm{sinc}(2\pi \tilde{x}).\] 
That concludes the proof of Proposition \ref{prop17P2p}.$\blacksquare$

\subsection{Proof of Theorem \ref{transprofP2p}}\label{proof1}

The refocused transversed profile is given by
\[\begin{split}
\frac{d}{\tilde{d}_2+\tilde{d}_1}&\frac{\lambda_{oc}^{1-\alpha_M}}{\theta}H^{\alpha_M}_{x_0}\Big(\omega_0,x_0+\frac{\lambda_{oc}}{\theta}\tilde{x},L\Big)=\frac{\lambda_{oc}}{2\theta}\sum_{j,l=1}^{N}\mathcal{T}^l_j(\omega_0,L)\phi_l(\omega_0,x_0)\\
&\phi_l(\omega_0,x_0+\lambda_{oc}\tilde{x}/\theta)\frac{A_j^2 d}{2}\Big[1-\cos\Big(\sigma_j\frac{2d_M+\lambda_{oc}^{\alpha_M}(\tilde{d}_2-\tilde{d}_1)}{d}\Big)\textrm{sinc}\Big(\sigma_j\frac{\lambda_{oc}^{\alpha_M}(\tilde{d}_2+\tilde{d}_1)}{d}\Big)\Big],
\end{split}\]
Using the following probabilistic representation
\[\mathcal{T}^l_j(\omega_0,L)=\E\Big[e^{-\int_0^L \Lambda^c_{Y^{N}_{s}}(\omega)\1_{(Y^{N}_{s}=N)\cup(Y^{N}_{s}=N-1)}ds}\1_{(Y^{N}_{L}=j)}\Big \vert Y^{N}_{0}=l\Big],\]
where  $\big(Y^{N}_{t}\big)_{t\geq0}$ is a jump Markov process, with state space $\{1,\dots,N\}$, intensity matrix $\Gamma^{c}(\omega_0)$, and invariant measure $\mu_{N}$, the uniform distribution over $\{1,\dots,N\}$, we have
\[\begin{split}
\Big\vert  \frac{\lambda_{oc}}{2\theta}\sum_{j,l=1}^{N}\mathcal{T}^l_j(\omega_0,L)&\phi_l(\omega_0,x_0)\phi_l(\omega_0,x_0+\lambda_{oc}\tilde{x}/\theta)\\
&\times \frac{A_j^2 d}{2}\cos\Big(\sigma_j\frac{2d_M+\lambda_{oc}^{\alpha_M}(\tilde{d}_2-\tilde{d}_1)}{d}\Big)\textrm{sinc}\Big(\sigma_j\frac{\lambda_{oc}^{\alpha_M}(\tilde{d}_2+\tilde{d}_1)}{d}\Big) \Big\vert\\ 
&\leq \lambda_{oc}^{1-\alpha_M}N \Big[\sum_{j=2}^{N}\frac{1}{\pi(j-1)}\mathbb{P}_{\mu_{N}}\big(Y^{N}_L=j\big)+\frac{1}{\sigma_1}\mathbb{P}_{\mu_{N}}\big(Y^{N}_L=1\big)\Big] \\
&\leq K \lambda_{oc}^{1-\alpha_M}\ln(N).
\end{split}\]
Consequently, the transverse profile of the refocused wave is given by 
\[\frac{\lambda_{oc}}{2\theta}\sum_{j,l=1}^{N}\frac{A_j^2 d}{2}\mathcal{T}^l_j(\omega_0,L)\phi_l(\omega_0,x_0)\phi_l(\omega_0,x_0+\lambda_{oc}\tilde{x}/\theta).\]
Let $\eta >0$ such that $\eta \ll1$. Thanks to \eqref{preveq}, we have 
\[\begin{split}
\frac{\lambda_{oc}}{2\theta}\sum_{j,l=1}^{N}\frac{A_j^2 d}{2}\mathcal{T}^l_j(\omega_0,L)&\phi_l(\omega_0,x_0)\phi_l(\omega_0,x_0+\lambda_{oc}\tilde{x}/\theta)\\
=&\frac{\lambda_{oc}}{2\theta}\sum_{j,l=1}^{[N(1-\eta)]}\frac{A_j^2 d}{2}\mathcal{T}^l_j(\omega_0,L)\phi_l(\omega_0,x_0)\phi_l(\omega_0,x_0+\lambda_{oc}\tilde{x}/\theta)\\
&+\mathcal{O}(\eta)\\
=&\frac{\lambda_{oc}}{2\theta}\sum_{j,l=1}^{[N(1-\eta)]}\mathcal{T}^l_j(\omega_0,L)\phi_l(\omega_0,x_0)\phi_l(\omega_0,x_0+\lambda_{oc}\tilde{x}/\theta)\\
&+\mathcal{O}(\eta).
\end{split}\]
Let $f^{\eta}(v)=\textbf{1}_{[0,1-\eta]}(v)$, we have
\[\begin{split}
\frac{\lambda_{oc}}{2\theta}\sum_{j,l=1}^{N}\frac{A_j^2 d}{2}\mathcal{T}^l_j(\omega_0,L)&\phi_l(\omega_0,x_0)\phi_l(\omega_0,x_0+\lambda_{oc}\tilde{x}/\theta)\\
=&\frac{\lambda_{oc}}{2\theta}\sum_{l=1}^{[N(1-\eta)]}\mathcal{T}^l_{f^{\eta}}(\omega_0,L)\phi_l(\omega_0,x_0)\phi_l(\omega_0,x_0+\lambda_{oc}\tilde{x}/\theta)\\
&+\mathcal{O}(\eta).
\end{split}\]
Now, we are able to use the continuous diffusion approximation given in \cite[Theorem 6.2]{gomez2}. In what follows $\mathcal{T}_{f^{\eta}}$ is the solution of the diffusion equation given in Theorem \ref{transprofP2p} but with initial conditions $f^{\eta}$. As a result, we have
\[\begin{split}
\frac{1}{N}\sum_{l=1}^{[N(1-\eta)]}\Big\lvert\mathcal{T}^l_{f^{\eta}}(\omega_0,L)&-\mathcal{T}_{f^{\eta}}\big(L,l/N\big) \Big\rvert\\
&\leq \sum_{l=1}^{[N(1-\eta)]-1}\int_{l/N}^{(l+1)/N}\big\lvert  \mathcal{T}^{[Nu]}_{f^{\eta}}(\omega_0,L)-\mathcal{T}_{f^{\eta}}\big(L,[Nu]/N\big) \big\rvert du\\
&\leq \int_0^{1}\big\lvert  \mathcal{T}^{N}_{f^{\eta}}(\omega_0,L,u)-\mathcal{T}_{f^{\eta}}(L,u)\big\vert du \\
&\quad+ \int_0^1\Big\lvert \mathcal{T}_{f^{\eta}}(L,u)- \mathcal{T}_{f^{\eta}}\big(L,[Nu]/N\big) \big\rvert du,
\end{split}\]   
where the terms on the right side of the last inequality converge to $0$ as $\omega_0\to+\infty$. Then
\[\begin{split}
\frac{\lambda_{oc}}{2\theta}\sum_{j,l=1}^{N}\frac{A_j^2 d}{2}\mathcal{T}^l_j(\omega_0,L)&\phi_l(\omega_0,x_0)\phi_l(\omega_0,x_0+\lambda_{oc}\tilde{x}/\theta)\\
=&\frac{\lambda_{oc}}{2\theta}\sum_{l=1}^{[N(1-\eta)]}\mathcal{T}_{f^{\eta}}(L,l/N)\phi_l(\omega_0,x_0)\phi_l(\omega_0,x_0+\lambda_{oc}\tilde{x}/\theta)\\
&+\mathcal{O}(\eta).
\end{split}\]
Moreover, we have 
\[\phi_j (\omega_0,x_0)\phi_j\big(\omega_0,x_0+\frac{\lambda_{oc}}{\theta}\tilde{x}\big)=\frac{A^2_j}{2}\Big[\cos\Big(\sigma_j\frac{\lambda_{oc}}{\theta d}\tilde{x}\Big)-\cos\Big(\sigma_j \frac{2x_0 +\lambda_{oc}\tilde{x}/\theta}{d}\Big)\Big],\]
and 
\[\begin{split}
\cos\Big(\sigma_j \frac{2x_0 +\lambda_{oc}\tilde{x}/\theta}{d}\Big)=&\cos\Big((\sigma_j-j\pi) \frac{2x_0 +\lambda_{oc}\tilde{x}/\theta}{d}\Big)\cos\Big(j\pi \frac{2x_0 +\lambda_{oc}\tilde{x}/\theta}{d}\Big)\\
&-\sin\Big((\sigma_j-j\pi) \frac{2x_0 +\lambda_{oc}\tilde{x}/\theta}{d}\Big)\sin\Big(j\pi \frac{2x_0 +\lambda_{oc}\tilde{x}/\theta}{d}\Big).
\end{split}\]
Using the Abel transform, Lemma \ref{coefgP2}, \eqref{preveq}, and the continuity of $v\mapsto\mathcal{T}_{f^{\eta_1}}(L,v)$ on $[0,1]$, we get
\[\lim_{\omega_0\to +\infty}\lambda_{oc}\Big\lvert \sum_{l=1}^{[N(1-\eta)]}\mathcal{T}_{f^{\eta}}(L,l/N)A^2_l\cos\Big(\sigma_j \frac{2x_0 +\lambda_{oc}\tilde{x}/\theta}{d}\Big) \Big\rvert=0.\]
Moreover, using \eqref{preveq} on more time and the fact that $\lim_{\omega_0} \sup_{j}\lambda_{oc}\lvert\sigma_j -j\pi \rvert=0$,
\[\begin{split}
\lim_{N\to +\infty}\frac{\lambda_{oc}}{2\theta}\sum_{l=1}^{[N(1-\eta)]}\mathcal{T}_{f^{\eta}}(L,l/N)&A^2_l\cos\Big(\sigma_j \frac{\lambda_{oc}}{\theta d}\tilde{x}\Big)\\
&=\lim_{N\to +\infty}\frac{\lambda_{oc}}{2\theta d^2}\sum_{l=1}^{[N(1-\eta)]}\mathcal{T}_{f^{\eta}}(L,l/N)\cos\Big(2\pi\frac{l}{N}\tilde{x}\Big)\\
&=(1-\eta)\int_0^{1-\eta}\mathcal{T}_{f^{\eta}}(L,u)\cos(2\pi u \tilde{x})du.
\end{split}\]  
Consequently, according to the proof of \cite[Theorem 6.2]{gomez2}, we have 
 \[\| \mathcal{T}_{f^{\eta}}(L,.)- \mathcal{T}_{1}(L,.)\|_{L^2([0,1])}\leq \| f^{\eta}- 1\|_{L^2([0,1])},\] 
and then
\[\begin{split}
\varlimsup_{\omega_0\to +\infty}\Big\lvert \frac{\lambda_{oc}}{2\theta}\sum_{j,l=1}^{N}\frac{A_j^2 d}{2}\mathcal{T}^l_j(\omega_0,L)\phi_l(\omega_0,x_0)&\phi_l(\omega_0,x_0+\lambda_{oc}\tilde{x}/\theta) \\
&-\int_0^1 \mathcal{T}_1(L,v)\cos(2\pi u \tilde{x})du \Big\rvert\leq K\eta.
\end{split}\]
This concludes the proof of Theorem \ref{transprofP2p}.$\blacksquare$

\subsection{Proof of Theorem \ref{transprof2P2p}}\label{proof2}

Let $M\geq 1$ and $f^M=\sum_{j=1}^M\alpha_j\phi_j\in\mathcal{E}_M$. Moreover, let
\[\forall x\in[0,d],\quad\tilde{f}^M(x)=\sum_{j=1}^M \alpha_j \phi_j(\omega_0,x).\]
Using \eqref{preveq} and 
\[\sup_{j\in\{1,\dots,M\}}\lvert\sigma_j -j\pi \rvert=\mathcal{O}\left(\frac{1}{N}\right), \] 
we have
\[\sup_{x\in[0,d]}\lvert f^M(x)-\tilde{f}^M(x)\rvert=\mathcal{O}\left(\frac{1}{N}\right).\]
Finally, by considering
\[\tilde{H}^1_{x_0}(\omega_0,x,L)=\frac{\tilde{d}_2+\tilde{d}_1}{d}\frac{\lambda_{oc}}{2\theta}\sum_{j,l=1}^{N}\frac{A_j^2 d}{2}\mathcal{T}^l_j(\omega_0,L)\phi_l(\omega_0,x_0)\phi_l(\omega_0,x),\]
we have
\[\begin{split}
\big<H^1_{x_0}(\omega_0,.,L)-\tilde{H}^1_{x_0}(\omega_0,.,L),f^M\big>_{L^2(0,d)}=&\big<H^1_{x_0}(\omega_0,.,L)-\tilde{H}^1_{x_0}(\omega_0,.,L),f^M-\tilde{f}^M\big>_{L^2(0,d)}\\
&+\big<H^1_{x_0}(\omega_0,.,L)-\tilde{H}^1_{x_0}(\omega_0,.,L),\tilde{f}^M\big>_{L^2(0,d)},
\end{split}\]
with
\[\begin{split}
\big\lvert \big<H^1_{x_0}(\omega_0,.,L)-\tilde{H}^1(\omega_0,.,L),f^M-&\tilde{f}^M\big>_{L^2(0,d)}\big\rvert \\ &\leq \frac{K}{N}N \Big[\sum_{j=2}^{N}\frac{1}{\pi(j-1)}\mathbb{P}_{\mu_{N}}\big(Y^{N}_L=j\big)+\frac{1}{\sigma_1}\mathbb{P}_{\mu_{N}}\big(Y^{N}_L=1\big)\Big]  \\
&\leq K\frac{\ln(N)}{N},
\end{split}\]
and
\[\begin{split}
\big\lvert \big<H^1_{x_0}(\omega_0,.,L)-&\tilde{H}^1_{x_0}(\omega_0,.,L),\tilde{f}^M\big>_{L^2(0,d)}\big\rvert \\
&\leq K \sum_{l=1}^M\Big[\mathcal{T}^l_{f^\eta}(\omega_0,L)+\frac{1}{[N\eta]+1} \mathbb{P}\big(Y^N_L\in\{[N\eta]+1,\dots,N\}\big\vert Y^N_0=l\big)\Big]
\end{split}\]
for $\eta>0$, and $f^\eta (v)=\textbf{1}_{[0,\eta]}$. Therefore, it suffices to study $\sum_{l=1}^M\mathcal{T}^l_{f^\eta}(\omega_0,L)$. To do this, let $g^{\eta}$ be a smooth function with compact support included in $[0,2\eta]$ and such that $0\leq f^{\eta}\leq g^{\eta}\leq f^{2\eta}$. Using the second part of \cite[Theorem 6.2]{gomez2},
\[\varlimsup_{N\to +\infty} \sum_{l=1}^M\mathcal{T}^l_{g^\eta}(\omega_0,L)=M\mathcal{T}_{g^\eta}(L,0)\leq M\overline{\mathbb{P}}_{0} \big( x(t)\in[0,2\eta]\big).\]
Here, we have used the probabilistic representation of $\mathcal{T}^l_{g^\eta}(\omega_0,L)$ introduced in the proof of \cite[Theorem 6.2]{gomez2}, where $\overline{\mathbb{P}}_{0}$ is the unique solution of a martingale problem starting from $0$. However, the probabilistic representation can be chosen such that the associated diffusion process has aabsolutely continuous transition probabilities with respect to the Lebesgue measure \cite{friedman}. Therefore,
\[\lim_{\omega_0\to +\infty}\big<H^1_{x_0}(\omega_0,.,L)-\tilde{H}^1_{x_0}(\omega_0,.,L),\tilde{f}^M\big>_{L^2(0,d)}=0,\]
and the rest of the proof is the same as the one of Theorem \ref{transprofP2p}.
$\blacksquare$

\subsection{Proof of Proposition \ref{propsinctrP2p}}\label{proof3}

Following the proof of Theorem \ref{transprofP2p} and using \cite[Theorem 6.4]{gomez2}, we get
\[\lim_{\omega_0 \to +\infty}\frac{\lambda_{oc}^{1-\alpha_M}}{\theta}H^{0,\alpha_M}_{x_0}\Big(\omega_0,x_0+\frac{\lambda_{oc}}{\theta}\tilde{x},L\Big)=\frac{\tilde{d}_2+\tilde{d}_1}{d}\int_0^1 \mathcal{T}_1(L,u)\cos(2\pi u\tilde{x})du,\]
where $\mathcal{T}_1(z,v)$ is a solution of 
\[\frac{\partial}{\partial z} \mathcal{T}_1 (z,u)= \frac{\partial}{\partial u}\left(a_{\infty}(\cdot)\frac{\partial}{\partial u}  \mathcal{T}_1\right)(z,u), \]
with boundary conditions 
\[ \frac{\partial}{\partial u} \mathcal{T}_1 (z,0)=0,\quad   \frac{\partial}{\partial u}\mathcal{T}_1 (z,1)=0, \quad \text{and}\quad \mathcal{T}_1 (0,u)=1.\]
However, thanks to the conservation energy caused by the two reflecting boundary conditions, and because the initial condition is constant equal to one, this problem admits only one solution, which is $\mathcal{T}_1(z,u)=1$. 

The transverse profile of the refocused wave is studied using the diffusive continuous model introduced in \cite{gomez2}, with two reflecting boundary conditions at $u=0$ (the top of the waveguide) and $u=1$ (the bottom of the waveguide). Here, the two reflecting boundary conditions mean that there is no radiative loss anymore (see Figure \ref{perterad2P2p}), and then the energy is conserved. This is the reason why $\mathcal{T}_1(z,u)=1$. Consequently, the sinc profile obtained in Proposition \ref{propsinctrP2p} is the best transverse profile that we can obtain.  $\blacksquare$

\end{document}